%% file: main.tex
\documentclass[11pt]{amsart}
\usepackage[T1]{fontenc}
\usepackage{mathptmx}
\usepackage{amsmath,amssymb,amsthm,amscd,mathtools}
\DeclareMathAlphabet{\mathcal}{OMS}{cmsy}{m}{n}
\usepackage{array,booktabs,longtable}
\usepackage{microtype}
\usepackage[hidelinks]{hyperref}
\usepackage{stackengine}
\stackMath
\usepackage{xcolor}

\input{macros}

\theoremstyle{plain}
\newtheorem{theorem}{Theorem}[section]
\newtheorem{proposition}[theorem]{Proposition}
\newtheorem{lemma}[theorem]{Lemma}
\newtheorem{corollary}[theorem]{Corollary}
\theoremstyle{definition}
\newtheorem{definition}[theorem]{Definition}

\theoremstyle{remark}
\newtheorem{remark}[theorem]{Remark}

\title{An equivariant index theorem for the Heisenberg calculus}
\author{Alexander Gorokhovsky}
\address{Department of Mathematics, University of Colorado, Boulder, Colorado 80309, USA}
\email{gorokhov@colorado.edu}
\author{Erik van Erp}
\address{Department of Mathematics, Dartmouth College, Hanover, New Hampshire 03755, USA}

\begin{document}
\begin{abstract}
We derive a local fixed-point formula for the equivariant index of Heisenberg elliptic operators on closed contact manifolds equipped with an action of a compact Lie group $G$ preserving the contact form. 
The formula expresses the character of the equivariant index in terms of a cyclic-cohomological character of the noncommutative Heisenberg principal symbol. 
The main new ingredient is the construction of traces $\Tr_\gamma$, for $\gamma\in G$, defined by integration over twisted conjugacy classes in the Heisenberg group. 
These generalize the trace used in our earlier nonequivariant index theorem \cite{GvE20}. The proof proceeds by a $K$-theoretic reduction to the equivariant index theorem for Toeplitz operators which we proved in \cite{GvE26}.
\end{abstract}
\maketitle

\setcounter{tocdepth}{1}
\tableofcontents

\section{Introduction}

On a closed contact manifold, a Heisenberg elliptic operator is a Fredholm operator.
In \cite{GvE20} we gave a local formula for the index.
The aim of this paper is to generalize this formula to the equivariant case where a compact Lie group acts, preserving the contact form.

\subsection{The index problem}
Let $M$ be a compact manifold without boundary of dimension $2n+1$, with contact form $\contact$ and contact hyperplane bundle $H=\ker\contact$.
Following Folland and Stein \cite{FS74}, the Heisenberg calculus \cite{Ta84,BG88} gives the Reeb field order $2$ and the vector fields in $H$ order $1$.
Principal symbols of order zero form a noncommutative algebra \[\Symb=\Psi^0_H(M)/\Psi^{-1}_H(M).\]
An operator with invertible principal symbol, a \emph{Heisenberg elliptic} operator, has a parametrix and is hypoelliptic. On a compact manifold it is Fredholm.
Motivating examples come from several complex variables. 
On the boundary of a strictly pseudoconvex domain the Kohn Laplacian $\Box_b$ on $(0,q)$-forms is Heisenberg elliptic for  $0<q<n$. The Szeg\"{o} projection and the Toeplitz operators of Boutet de Monvel \cite{Bo79} are operators of order zero in the calculus.

The index problem for these operators was studied extensively by Epstein and Melrose \cite{EM98,Ep04,EMxx}.
Their formula for Hermite operators generalizes Boutet de Monvel's index formula for Toeplitz operators \cite{Bo79}, but their work found obstacles in the way of a general formula (see the discussion in \cite{Ep04}).
$K$-theoretic index theorems were obtained in \cite{vE10a, vE10b}.
These results do not readily yield formulas computable from local data, because the relevant isomorphisms in analytic $K$-theory are not given by explicit formulas.
A step forward was the explicit calculation in \cite{BvE14} of the class in geometric $K$-homology determined by a Heisenberg elliptic operator. However, the construction of this $K$-cycle involves an initial deformation of the symbol which is not canonical. This obstructs its generalization to the equivariant case.

In \cite{GvE20} we gave a local formula for general Heisenberg elliptic operators,
\begin{equation}\label{eq:non-equivariant index formula}
     \ind P=\int_M \chi(\sigma_H(P))\wedge \Ahat(M)
\end{equation}
The Heisenberg symbol of $P$ determines a class in $K$-theory,
\[ [\sigma_H(P)]\in K_1(\Symb)\]
The homomorphism
\[\chi:K_1(\Symb)\to H^{\mathrm{odd}}(M;\CC)\]
is the character of a generalized cycle in cyclic cohomology \cite{gor1,gor2}, given by an explicit expression  in the  symbol, a connection on $H$ and its curvature (see (\ref{eq:chi non-equivariant}) below).
A crucial ingredient of the formula is a regularized trace on $\Symb$,
for which we subsequently found a much simpler description in \cite{GvE25}.

For elliptic operators, the step after the Atiyah--Singer index theorem was the equivariant index theorem and the Lefschetz fixed-point formula of Atiyah, Segal and Singer \cite{AS2,AS3}.
This paper answers the corresponding question for the Heisenberg calculus.

If $G$ acts on $M$ preserving $\contact$, the index of a $G$-equivariant Heisenberg elliptic operator $P$ is an element in the representation ring of $G$,
\[ \Ind_G (P)\in R(G)\]
Our formula for the character of the equivariant index is local in the principal symbol and localizes on the fixed-point set $M^\gamma$.
For Toeplitz operators, we proved an equivariant index formula  in \cite{GvE26}. This formula is the input for our proof  here.

\subsection{The main theorem}
Throughout, $G$ is a compact Lie group acting on $M$ and preserving $\contact$.
We fix a $G$-invariant complex structure $J$ on $H$ compatible with the symplectic form $\omega=d\contact|_H$.
This  makes $H$ a Hermitian vector bundle which we denote by $\Hzo$.
The fixed-point set $M^\gamma$ is a contact submanifold.
Along $M^\gamma$ we have $\Hzo|_{M^\gamma}=F\oplus N^\gamma$, where $F=\ker(1-\gamma)$ is the complex contact bundle of $M^\gamma$ and its orthogonal complement $N^\gamma$ is the normal bundle of $M^\gamma$ in $M$.
$N^\gamma$ is a complex vector bundle on which $\gamma$ acts unitarily without eigenvalue $1$.

The heart of the formula is an equivariant  character of the principal symbol,
\[ \chi_\gamma:K_1^G(\Symb)\to H^{\mathrm{odd}}(M^\gamma;\CC)\]
Like $\chi$, it is given by an explicit local expression in the symbol, a $G$-invariant connection on $\Hzo$ and its curvature, now with values in differential forms on $M^\gamma$.

\begin{theorem}\label{thm:main-intro}
Let $P$ be a $G$-equivariant Heisenberg elliptic operator of order zero acting on sections of a
$G$-equivariant Hermitian vector bundle.
For every $\gamma\in G$,
\[
 \Tr\bigl(\gamma\,|\,\Ind_G(P)\bigr)
   =\int_{M^\gamma}
      \chi_\gamma\bigl(\sigma_H(P)\bigr)
      \wedge \Ahat(M^\gamma)\wedge \frac{1}{\Delta_\gamma(N^\gamma)}
\]
where each connected component of $M^\gamma$, of dimension $2r+1$, is oriented by $\contact\wedge(d\contact)^r>0$.
\end{theorem}

$\Ahat(M^\gamma)$ is the $\widehat A$-class of the tangent bundle of $M^\gamma$.
The characteristic class $\Delta_\gamma(N^\gamma)$ of the complex vector bundle $N^\gamma$ is analogous to the denominator that appears in Atiyah-Segal-Singer fixed-point formulas.

\begin{definition}\label{def:Delta_gamma}
Let $N$ be  a complex vector bundle with a unitary automorphism $\gamma$ with locally constant eigenvalues and multiplicities, and without eigenvalue $1$.
Denote the eigenvalues by $e^{i\vartheta_j}$ with $0<\vartheta_j<2\pi$, and let $x_j$ be the   Chern roots on the respective eigenbundles. Then 
\begin{equation*}
 \Delta_\gamma(N)
 :=
 \prod_{j}
   \bigl(e^{i\vartheta_j/2}e^{x_j/2}-e^{-i\vartheta_j/2}e^{-x_j/2} \bigr)
\end{equation*}
\end{definition}

\subsection{Orbital traces}\label{sec:intro traces}
The new idea of the paper is the construction of the traces appearing in the character $\chi_\gamma$.

In \cite{GvE20} a symbol is a pair $\sigma=(\sigma_+,\sigma_-)$ of elements of Weyl algebras of the fibers of $H^*$. The trace is 
\[\tau(\sigma)=\Trh(\sigma_+)+(-1)^{n+1}\Trh(\sigma_-)\]
where $\Trh$ is the zeta-regularized trace of Epstein and Melrose \cite{EMxx}.
In \cite{GvE25} we found a simpler description.
By Fourier transform, an order zero symbol at $x\in M$ is a regular distribution (i.e. smooth off the identity) on the Heisenberg group 
\[ \Heis_x=H_x^{1,0}\times \RR\]
homogeneous of degree $-2n-2$ for the parabolic dilations.
The product of symbols is convolution.
In this kernel picture, evaluation at a central element $(0,c)$, $c\ne0$, is a trace $\orbtr{c}(f)=f(0,c)$, and the trace $\tau$ is
\[ \tau = \frac{(2\pi i)^{n+1}}{2(n!)}\Bigl(\orbtr{+1}-(-1)^n\orbtr{-1}\Bigr)\]
This simplified picture of the trace in \cite{GvE25} is what makes the equivariant problem tractable.

At a point $x\in M^\gamma$ in the fixed-point set, $\gamma$ acts on the Heisenberg group $\Heis_x$ by a group automorphism that fixes the center.
We replace the central element $(0,c)$ by its orbit $\Gamma_c$ under $\gamma$-twisted conjugation $h\mapsto a\,h\,\gamma(a^{-1})$.
Integration over conjugacy classes  gives traces.
The conjugacy class of a central element consists of a single point.
Integration over twisted conjugacy classes gives twisted traces.
Restriction of a twisted trace to invariant elements gives a trace.

Concretely, the orbit $\Gamma_c$ is the graph of $z\mapsto c+q(z)$ over the normal space $N^\gamma_x$, with $q$ a quadratic form determined by $\gamma$ (Section~\ref{sec:twisted orbital traces on Heis}).
Integration over the orbit $\Gamma_c$ is a $\gamma$-twisted trace,
\[ \orbtr{c}(f)=\int_{N^\gamma_x}f(z,c+q(z))\,dz\qquad \orbtr{c}(a*b)=\orbtr{c}(b*a^\gamma)\]
In particular $\orbtr{c}$ is a trace on $G$-invariant symbols.
The construction is defined fiberwise over $M^\gamma$, and therefore gives a trace  with values in functions on $M^\gamma$,
\[ \Tr_\gamma:\Symb(M)^G\to C^\infty(M^\gamma)\qquad \Tr_\gamma = \kappa_\gamma \cdot \frac{(2\pi i)^{r+1}}{2(r!)}\Bigl(\orbtr{+1}-(-1)^r\orbtr{-1}\Bigr):\]
Here $\kappa_\gamma$ is an explicit constant determined by the eigenvalues of $\gamma$ on $N^\gamma$ (Definition~\ref{def:Tr_gamma}).
The localization on $M^\gamma$ is built into the construction.
If $\gamma$ acts trivially, the orbits $\Gamma_c$ are the points $(0,c)$ and $\Tr_\gamma$ agrees with $\tau$.

\subsection{The character}
The  character $\chi_\gamma$ is constructed as in \cite{GvE20}, replacing the trace $\tau$ by $\Tr_\gamma$.
 Choose a $G$-invariant unitary connection $\nabla$ on $H^{1,0}$. It determines a connection on the symbol algebra $\Symb$ that preserves invariant forms.
We enlarge the symbol algebra,
\[\Symb\subset \Alg_H\subset \Omega^\bullet(\Alg_H)=\Alg_H\otimes_{C^\infty(M)} \Omega^\bullet(M)\]
so that curvature is an inner derivation,
\[\con^2(a)=[\curv,a]\qquad \curv\in\Omega^2(\Alg_H)\]
This provides the curved differential graded algebra needed for the construction of a cycle in cyclic homology.
In the kernel picture the enlarged algebra is obtained by formally inverting the central distribution $\Tsf=-\partial_t\delta_e$.
The traces $\Tr_\gamma$ extend to $G$-invariant $\Alg_H(M)$-valued forms, with values in $\Omega^\bullet(M^\gamma)$, and 
\[\Tr_\gamma(\con\eta)=d\,\Tr_\gamma(\eta)\]
The closed graded trace  and curved differential graded algebra give a morphism from cyclic chains to differential forms. Applying this morphism to the Connes-Chern character $\ch(\sigma)$ gives the character $\chi_\gamma(\sigma)$, 
\begin{multline}\label{eq:chi non-equivariant}
\chi_\gamma(\sigma) = \sum \limits_{l\ge 0} \sum \limits_{i_0, \ldots, i_{2l+1} \ge 0} \left(\frac{-1}{2 \pi i}\right)^{I+l+1} \frac{l!}{(I+2l+1)!}\\
 \Tr_\gamma \left(  \sigma^{-1} \curv^{i_0} \con(\sigma)\curv^{i_1}\con(\sigma^{-1}) \dots \con(\sigma) \curv^{i_{2l+1}}\right)
\end{multline}
\[I=i_0+i_1+\cdots+i_{2l+1}\]
If $\gamma$ acts trivially (e.g.\ if $\gamma=e$) then $\Tr_\gamma=\tau$ is the trace of  \cite{GvE25}, and $\chi_e=\chi$ is the non-equivariant character of \cite{GvE20}.

\subsection{The proof}
The proof follows the $K$-theoretic reduction of \cite{GvE20}.
Rather than proving Theorem \ref{thm:main-intro} directly for arbitrary Heisenberg elliptic operators, we use $K$-theory to reduce the problem to Toeplitz operators.
We show that $K_1^G(\Symb)$ is generated by symbols of multiplication operators and Toeplitz operators.
Both sides of the index formula vanish on multiplication operators.
Thus, it suffices to verify that our formula is correct for Toeplitz operators.

Let $T_f$ be the Toeplitz operator corresponding to a $G$-equivariant vector bundle automorphism $f$.
The direct sum  $\tilde T_f=T_f\oplus 1$  with the identity on the orthogonal complement of the range of the Szeg\"{o} projection is a $G$-equivariant order zero Heisenberg elliptic operator with the same index as $T_f$.
In Section \ref{sec:Toeplitz operators} we calculate the character $\chi_\gamma(\sigma_H(\tilde T_f))$ and show that in $H^\bullet(M^\gamma;\CC)$ we have equality
\[
 \chi_\gamma(\sigma_H(\tilde T_f))\wedge \Ahat(M^\gamma)\wedge \frac{1}{\Delta_\gamma(N^\gamma)}
 =\ch_\gamma(f)\wedge
   \Td(M^\gamma)\wedge \frac{1}{\ch_\gamma(\lambda_{-1}N^*)}
\]
The left-hand side is the integrand in Theorem~\ref{thm:main-intro}, while the right-hand side is the integrand in the equivariant index formula for Toeplitz operators  which we proved in \cite{GvE26}.
Thus Theorem \ref{thm:main-intro} is correct when restricted to Toeplitz operators, and therefore holds for all Heisenberg elliptic operators of order zero.

\subsection{Outline of the paper}
Section~\ref{sec:symbols} reviews the kernel picture of Heisenberg symbols.
Section~\ref{sec:enlarged-symbol-algebra} constructs the enlarged  symbol algebra $\Symb\subset \Alg_H$ and its curvature $\curv\in\Omega^2(\Alg_H)$.
Section~\ref{sec:orbital-traces} constructs the twisted orbital traces from which $\Tr_\gamma$ is built.
Section \ref{sec:character map} defines  the character map $\chi_\gamma$. 
Section \ref{sec:Toeplitz operators} computes the character $\chi_\gamma(\sigma_H(T_f))$ of the Heisenberg symbol of a Toeplitz operator $T_f$, and verifies that our general formula here agrees with the equivariant index theorem of \cite{GvE26}.
Our proof of the index formula is in Section~\ref{sec:index-proof}.
The Appendix collects our conventions and compares them with those of \cite{GvE20,GvE25,GvE26}.

\section{Symbols as distributions}\label{sec:symbols}
In \cite{GvE20} we represented principal symbols in the Heisenberg calculus as pairs of elements in the Weyl algebra.
In this paper we use the dual picture, representing symbols as distributions on the Heisenberg group with convolution product.
The two pictures yield isomorphic symbol algebras via Fourier transform.
This section briefly reviews conventions and fixes notation.

For details on the Heisenberg pseudodifferential calculus see \cite{Ta84,BG88, CGGP92, Ep04, Po08}.
A brief introduction is in \cite[Sections~2--3]{GvE20}.
Our approach here follows the treatment in \cite{vEY19}.

\subsection{The Heisenberg group}\label{sec:Heisenberg group}
Let $V$ be a complex vector space of complex dimension $n$ with Hermitian form
\[\langle u, v \rangle =g(u, v)-i \omega (u, v)\]
Here $g$ is a positive definite real inner product, and $\omega$ is a symplectic form.
We assume that the hermitian form is linear in the first variable, and
\[g(Ju,v)=\omega(u,v)\qquad Ju=iu\]
Denote by $\Heis$ the Heisenberg group $V \times \RR$ with product
\begin{equation}\label{eq:heisenberg group law}
(z, t) (w, s)= (z+w, t+s-\frac{1}{2}\omega(z, w)) =(z+w, t+s+\frac{1}{2}\operatorname{Im}\ang{z,w})
\end{equation}
The Haar measure is $\dd h=\dd z\,\dd t$, with $\dd z$ normalized as
\[
       \dd z=\frac{1}{n!}|\omega^n|
\]
The parabolic dilations are
\begin{equation*}
     \delta_a(z,t)=(az,a^2t)\qquad a>0
\end{equation*}
The homogeneous dimension of $\Heis$ is $Q=2n+2$.
Whenever we refer to homogeneous distributions or functions on $\Heis$ of degree $j$ we mean homogeneity with respect to the dilations $\delta_a$.
A function $f$ is homogeneous of degree $j$ if
\[ f(az,a^2t) = a^{j}f(z,t)\qquad a>0 \]
A distribution $f$
is homogeneous of degree $j$ if for every test function $\varphi\in C_c^\infty(\Heis)$,
\[ \ang{f,\varphi\circ \delta_a} = a^{-Q-j}\ang{f,\varphi}\]
The Dirac delta is homogeneous of degree $-Q$.

\subsection{Convolution of homogeneous distributions}\label{sec:model symbols}


We denote by $\Symb(\Heis)$ the set of regular distributions on $\Heis$ (i.e.\ smooth when restricted to $\Heis\setminus \{e\}$) that are homogeneous of degree $-Q=-2n-2$.

Convolution of two homogeneous distributions $f,g\in \Symb(\Heis)$ is defined as follows.
Let $\phi \in C_c^\infty\left(\Heis\right)$ be such that $\phi =1 $ in a neighborhood of $0$.
Then $f * g$ is, by definition, the principal part (of homogeneous degree $-Q$) in the asymptotic expansion of $\phi f * \phi g$. The result does not depend on the choice of $\phi$.

The Dirac delta is the unit of the algebra $\Symb(\Heis)$.

\subsection{Contact manifolds}\label{sec:contact}

Throughout this paper, $M$ is a compact smooth manifold without boundary of dimension $2n+1$,
equipped with a contact form $\contact$.
A compact Lie group $G$ acts on $M$, preserving $\contact$.

The Reeb field $\Reeb$ is the vector field on $M$ determined by $d\contact(\Reeb,\,-\,) = 0$, $\contact(\Reeb)=1$. The action of $G$ preserves $\Reeb$.

Set $H=\ker\contact \subset TM$ and
\[\omega:=d\contact|_H\]
(for the sign conventions of \cite{GvE20,GvE25,GvE26} see Appendix~\ref{sec:conventions}).
We orient $M$ by $\contact\wedge (d\contact)^n>0$.

We fix a $G$-invariant complex structure $J\in \mathrm{End}(H)$, $J^2=-\mathrm{Id}$,
compatible with the symplectic structure, i.e.\ $\omega(Jv,Jw) = \omega(v,w)$ and
\begin{equation*}
g(v,w) = \omega(v, Jw)
\end{equation*}
is positive definite.
Define a Hermitian form, linear in the first variable, by
\begin{equation*}
\ang{v,w}:= g(v,w)-i\omega(v,w)
\end{equation*}
We write $H^{1,0}$ for the resulting complex vector bundle with its Hermitian structure and $H^{0,1}$ for its conjugate.
The action of $G$ on $H$ preserves $J$, $g$ and $\omega$.
We extend $g$ to a $G$-invariant Riemannian metric on $M$ by $\Reeb\perp H$ and $|\Reeb|=1$.

Fix a point $x\in M$. We denote by $\Heis_x$ the Heisenberg group \[\Heis_x=H_x\times \RR\]
with group law (\ref{eq:heisenberg group law}) for the complex vector space $V=H_x^{1,0}$ with symplectic form $\omega=d\contact|_H$.
The bundle of Heisenberg groups is denoted
\[ T_HM=\bigsqcup_{x\in M} \Heis_x\]
As a fiber bundle over $M$, $T_HM$ identifies with the tangent bundle $TM=H\times \RR$.

For $\gamma\in G$, the restriction of $\contact$ to the fixed set $M^\gamma=\{x\in M\mid \gamma x=x\}$ is again a contact form.
The Reeb field $\Reeb$ is tangent to $M^\gamma$, and the normal bundle $N^\gamma$ of $M^\gamma$ in $M$ is a complex sub-bundle of $H^{1,0}|_{M^\gamma}$.
The action of $\gamma$ on $N^\gamma$ is unitary.

On each connected component of $M^\gamma$, the normal bundle decomposes into complex eigenbundles
\[ N^\gamma = \bigoplus_{\lambda\ne 1} N^\gamma_\lambda\qquad N_\lambda^\gamma=\ker (\gamma - \lambda)\]
These complex eigenbundles are mutually orthogonal for the hermitian metric, as well as symplectically orthogonal.

\subsection{Heisenberg symbols}\label{sec:Hsymb}
The Heisenberg calculus on a contact manifold is a $\ZZ$-filtered algebra $\Psi^\bullet_H(M)$ of pseudodifferential operators.
We denote the noncommutative algebra of Heisenberg principal symbols of order zero by $\Symb(M)$,
\[ 0\to \Psi_H^{-1}(M)\to \Psi_H^0(M)\to  \Symb(M)\to 0\]
An element in $\Symb(M)$ is a smooth section in a bundle whose fiber, at each point $x\in M$,
is the algebra $\Symb(\Heis_x)$ of section \ref{sec:model symbols}, with $\Heis_x=H^{1,0}_x\times \RR$.

This description differs from the Weyl pairs used in \cite{GvE20}.
The Fourier transform of a distributions $f\in\Symb(\Heis_x)$
is a smooth function $\sigma=\hat f$ on $H^*_x\times \RR\setminus \{(0,0)\}$ that is homogeneous of degree zero,
\[ \sigma(a\xi,a^2\ft) = \sigma(\xi,\ft)\qquad a>0\]
Then $\sigma$ is determined by its restrictions $\sigma_+,\sigma_-$
to the two hyperplanes $H^*_x\times \{-1\}$ and $H_x^*\times \{1\}$.
(For the labels used in \cite{GvE25} see Appendix~\ref{sec:conventions}.)
In this picture, multiplication of symbols is
\[ (\sigma_+,\sigma_-)(\rho_+,\rho_-)=(\sigma_+\#\rho_+,\rho_-\#\sigma_-)\]
where $\#$ is the Weyl algebra product of functions on the symplectic vector space $H^*_x$ (with symplectic form dual to $\omega=d\contact|_H$).

In this paper we do not represent Heisenberg symbols as Weyl pairs, but as distributions on the Heisenberg group.
The two algebras are isomorphic via Fourier transform.
The distribution picture significantly simplifies the construction of the trace that appears in our index formula.

\section{A curved dga}\label{sec:enlarged-symbol-algebra}

As in \cite{GvE20}, we need to enlarge the order-zero symbol algebra $\Symb(M)$ so that curvature is an inner derivation. The construction in this section is closely related to the construction in \cite[Sections~3--5]{GvE20},
except that here we work with distributions on the Heisenberg group rather than with Weyl pairs.
We suspect that the two constructions yield isomorphic algebras, but do not prove this, as none of our results depend on it.

\subsection{An enlarged algebra}

As above, let $\Heis=V\times \RR$.
For $j\geq0$, let $\cR^j(\Heis)$ be the space of regular distributions $f$
on $\Heis=V\times \RR$ that are homogeneous of degree $-Q-j$.
In particular,
\[ \Symb(\Heis)=\cR^0(\Heis)\]
Denote the direct sum as
\[ \cR(\Heis):=\bigoplus_{j\ge 0} \cR^{j}(\Heis)\]
It is a graded algebra under convolution.
The direct sum of the even degrees is denoted
\[ \cR_q(\Heis):=\bigoplus_{j\ge 0} \cR^{2j}(\Heis)\]
Define a map
\[
\heis\to \cR(\Heis)\qquad
 X\mapsto \mathsf X:=-X\delta_e
 \qquad
 \langle \mathsf X,\varphi\rangle=(X\varphi)(e)
\]
where $X\varphi$ is defined by representing $X$ as a translation invariant vector field on $\Heis$.
Evaluation at the identity $X\varphi(e)$ is the same for left and right-invariant vector fields.
In particular, we have
\[\Tsf=-\partial_t\delta_e\qquad \ang{\Tsf,\varphi} = \partial_t\varphi(e)\]
$\Tsf$ is regular and homogeneous of order $-Q-2$, and so $\Tsf\in \cR^2(\Heis)$.
$\Tsf$ is central in the convolution algebra.

The localization $\cR_q(\Heis)[\Tsf^{-1}]$ of $\cR_q(\Heis)$ at $\Tsf$ is a graded algebra,
with
\[\deg \Tsf^{-1}=-2\]
Algebraically we have a balanced tensor product,
\begin{equation}\label{eq:balanced tensor product}
\cR_q(\Heis)[\Tsf^{-1}]\cong \cR_q(\Heis)\otimes_{\CC[\Tsf]} \CC[\Tsf, \Tsf^{-1}]
\end{equation}
where the action of $\CC[\Tsf]$ on $\cR_q(\Heis)$ is determined by
\[ \Tsf f=-\partial_t f\]
Thus, the defining relation of the tensor product is
\[ -\partial_tf\otimes \Tsf^k = f\otimes \Tsf^{k+1}\]

\begin{definition}\label{eq:degree-zero-localized-algebra}
The algebra $\Alg_H(\Heis)$ is the degree-zero part of $\cR_q(\Heis)[\Tsf^{-1}]$,
\begin{equation*}
 \Alg_H(\Heis)
   :=\left(\cR_q(\Heis)[\Tsf^{-1}]\right)^0
\end{equation*}
\end{definition}
Elements of $\Alg_H(\Heis)$ are finite sums of terms $\Tsf^{-j}f_j$ with
$f_j\in\cR^{2j}(\Heis)$, subject to the relation
\[ \Tsf^{-j-1}(-\partial_tf)=\Tsf^{-j}f\]
The inverse $\Tsf^{-1}$ is formal.
We do not represent $\Tsf^{-1}$ as a distribution on $\Heis$.
Note that $\Tsf$ is not a zero-divisor of $\cR_q(\Heis)$, since $\Tsf * f=-\partial_t f=0$ implies $f$ is not in $\cR^j(\Heis)$ for any $j\ge 0$. Therefore the inclusion
\[\Symb(\Heis)\to \Alg_H(\Heis)\]
is injective.

\begin{remark}
The Fourier transform of $\Tsf$ is $-i\ft$, where $\ft$ is dual to the central variable. Localization at $\Tsf$ amounts to removing the $\ft=0$ plane.
In the Weyl pair picture of Heisenberg symbols, a term $\Tsf^{-j}f_j$ corresponds to a Weyl pair $(\sigma_+,\sigma_-)$ where $\sigma_\pm$ are asymptotically of degree $2j$ at infinity. Sums of such terms correspond to elements of the enlarged algebra $\Alg$ defined in \cite[Section 3.3]{GvE20}, which is the analogue of $\Alg_H(\Heis)$ in the Weyl pair picture.
\end{remark}

\subsection{Action of the symplectic Lie algebra}

The symplectic group $\mathrm{Sp}(V)$ acts by automorphisms on $\Heis=V\times \RR$ (acting trivially on $\RR$).
This determines an action on $\Alg_H(\Heis)$ by automorphisms.
Differentiation of this action gives an action of the Lie algebra $\mathfrak{sp}(V)$ on $\Alg_H(\Heis)$ by derivations.
In this section we show that $\mathfrak{sp}(V)$ acts by inner derivations.
See also \cite[Section 4]{GvE20}.

Let $\heis=V\oplus \RR$ be the Lie algebra of $\Heis=V\times \RR$ with
\[ [u,v] = -\omega(u,v)\,T\qquad u,v\in V\]
The Lie algebra $\heis$ is graded by $\heis_1=V$, $\heis_2=\RR$.
The universal enveloping algebra $\mcU(\heis)$ is graded accordingly.
Denote the central element by $T=(0,1)\in \heis$.
Let $S^pV$ denote the space of symmetric sums of monomials $v_1v_2\cdots v_p$ with $v_j\in V$.
We can symmetrize any product of elements in $V$ at the expense of factors of $T$, and get
\[ \mcU^p(\heis) = S^pV\oplus T\cdot S^{p-2}V\oplus T^2\cdot S^{p-4}V\oplus \cdots\]
We have
\[ [S^2V,S^2V]\subseteq T\cdot S^2V\qquad [S^2V,V]\subseteq T\cdot V\]
Note that $A\mapsto TA$ is an injective map $\mcU(\heis)\to \mcU(\heis)$.
Localizing at $T$ we get
\[\mcU(\heis)[T^{-1}]=\mcU(\heis)\otimes_{\CC[T]} \CC[T^{-1}, T] =\bigoplus_{p\ge 0, \;q\in\ZZ} T^q\cdot S^pV\]
In this algebra, consider the subspace
\[\Lg:=T^{-1}\cdot S^2V\]
Then $\Lg$ is a Lie algebra that acts on $V$ by commutators,
\[ [\Lg,\Lg]=\Lg\qquad [\Lg, V]\subset V \]

\begin{lemma}\label{lem:sp(V) action on V}
We have an isomorphism of Lie algebras
\[ \mu:\mathfrak{g}\cong \mathfrak{sp}(V)\qquad \mu(q) =  [q,\,\cdot\,] \]
\end{lemma}
\begin{proof}
Let $q\in \Lg$.
For any $u,v\in V$ the commutator $[u,v]\in \RR\cdot T$ is central, and so
\[ 0 = \mu(q)[u,v]=[\mu(q)u,v]+[u,\mu(q)v]\]
Then
\[ \omega(\mu(q)u,v)+\omega(u,\mu(q)v)=0\]
means that $\mu(q)\in \mathfrak{sp}(V)$.
The map $q\mapsto \mu(q)$ is injective, since no element in $S^2V$ commutes with all of $V$.
It is an isomorphism by $\dim \mathfrak{g} = \dim S^2V=\dim \mathfrak{sp}(V)$.
\end{proof}

The map
\[\heis\to \mcE'(\Heis)\qquad
 X\mapsto \mathsf X:=-X\delta_e\]
maps Lie brackets in $\heis$ to
commutators in the convolution algebra $\mcE'(\Heis)$.
Therefore the map extends to an algebra homomorphism
\[ \mcU(\heis)\to \mcE'(\Heis)\]
Denote the composition
\[ \mathfrak{sp}(V)\cong T^{-1}\cdot S^2V\subset \mcU(\heis)[T^{-1}]\to \mcE'(\Heis)[\Tsf^{-1}]\]
by
\[ \nu:\mathfrak{sp}(V)\to \mcE'(\Heis)[\Tsf^{-1}]\]
The symplectic group $\mathrm{Sp}(V)$ acts by automorphisms on the convolution algebra $\mcE'(\Heis)$.
The Lie algebra $\mathfrak{sp}(V)$ acts on $\mcE'(\Heis)$ by derivations.
Note that $A.\Tsf=0$ for $A\in\mathfrak{sp}(V)$. The action extends to $\mcE'(\Heis)[\Tsf^{-1}]$ by $A.\Tsf^{-1}=0$.

\begin{lemma}\label{lem:augmentation}
The augmentation
\[ \epsilon:\mcE'(\Heis)\to \CC\qquad \epsilon(f)=\ang{f,1}  \]
is a unital algebra homomorphism.
The kernel of $\epsilon$ is the ideal generated by $\heis$,
\[ \ker \epsilon = \heis * \mcE'(\Heis) = \{\sum_j \mathsf X_j*g_j\mid X_j\in \heis,\,g_j\in \mcE'(\Heis)\}\]
\end{lemma}
\begin{proof}
Recall the definition of convolution of distributions
\[\ang{f\ast g,\varphi}=\ang{f\otimes g, \varphi\circ m}\]
where $m:\Heis\times\Heis\to \Heis$ is multiplication in the group and
\[f\otimes g\in  \mcE'(\Heis\times\Heis)\]
This gives
\[ \epsilon(f*g)=\ang{f * g,1}=\ang{f\otimes g, 1\otimes 1}=\ang{f,1}\ang{g,1}=\epsilon(f)\epsilon(g)\]
Clearly $\epsilon(\delta_e)=1$. This proves that $\epsilon$ is unital algebra homomorphism.

For the kernel, the inclusion $\supseteq$ follows from $\epsilon(\mathsf X\ast g)=\epsilon(\mathsf X)\epsilon(g)$
and $\epsilon(\mathsf X)=\ang{\mathsf X,1}=(X\,1)(e)=0$.
To complete the proof we need to show $\subseteq$.

Define
\[ F_t:\Heis\to \Heis\qquad F_t(h)=\exp(t\log h)\]
Choose a basis $X_1,\dots, X_{2n}$ for $V$ and write $X_0=T$, and denote the corresponding exponential coordinates by
\[ x_j:\Heis\to \RR\qquad \log h=\sum_j x_j(h)X_j\]
Now fix $f\in \mcE'(\Heis)$.
For a test function $\varphi$, integration along the curve $t\mapsto F_t(h)$ gives
\[ \varphi(h) = \varphi(e)+\sum_j x_j(h) \int_0^1 (X_j\varphi)(F_t(h))\,dt\]
where $X_j$ is interpreted as a right-invariant vector field on $\Heis$.
Assuming $\epsilon(f)=0$ we get
\[
\ang{f,\varphi}
=\sum_j \ang{f, x_j \int_0^1 (X_j \varphi)\circ F_t \,dt}
\]
Set
\[ \ang{g_j,\varphi}:=\int_0^1 \ang{f, x_j\cdot (\varphi\circ F_t)}\,dt\]
Each $g_j$ is a well-defined distribution since $(t,h)\to x_j(h)\varphi(F_t(h))$ is smooth,
and
\[ \supp g_j\subset \bigcup_{t\in [0,1]} F_t(\supp f)\]
where the right-hand side is compact as the continuous image of $[0,1]\times \supp f$.

Now
\[
\ang{f,\varphi}
= \sum_j \ang{g_j, X_j \varphi}
= \sum_j \ang{\mathsf X_j*g_j, \varphi}
\]
So
\[f=\sum_j \mathsf X_j*g_j\in \heis * \mcE'(\Heis)\]
\end{proof}

\begin{lemma}\label{lem:vanishing D}
If a derivation $D$ of $\mcE'(\Heis)$ vanishes on $\heis$ then $D=0$.
\end{lemma}
\begin{proof}
Fix $f\in \mcE'(\Heis)$. Since $f-\epsilon(f)\delta_e$ is in the kernel of $\epsilon$, Lemma \ref{lem:augmentation} gives
\[ f - \epsilon(f)\delta_e =\sum_j \mathsf X_j * g_j\]
For any derivation we have $D(\delta_e)=0$. Since $D$ vanishes on $\heis$ we have $D(\mathsf X_j)=0$.
We get
\[ D(f) = \sum_j \mathsf X_j * D(g_j)\]
or
\[ \ang{D(f), \varphi} = \sum \ang{\mathsf X_j\ast D(g_j), \varphi}= \sum \ang{D(g_j), X_j \varphi}\]
We have $\ang{D(f), 1}=0$.
Induction on the degree of $p$ (for the dilations $\delta_a$) proves that $D(f)$ annihilates all polynomials on $\Heis$.
The values of $\ang{D(f),p}$ for polynomials $p$ are the Taylor coefficients at the origin
of the Fourier transform of $D(f)$.
This transform is entire because $D(f)$ is compactly supported.
Since $\ang{D(f),p}=0$ for all polynomials $p$ we get $D(f)=0$.
\end{proof}

\begin{proposition}\label{prop:sp(V) action is inner}
The action of $\mathfrak{sp}(V)$ on $\mcE'(\Heis)[\Tsf^{-1}]$ is inner,
\[ A.f=[\nu(A),f]\qquad A\in \mathfrak{sp}(V),\;f\in\mcE'(\Heis)[\Tsf^{-1}]\]
\end{proposition}
\begin{proof}
A priori, we have two derivations $f\mapsto A.f$ and $f\mapsto [\nu(A),f]$.
Since $\Tsf$ is central, both derivations annihilate $\Tsf$, and it suffices to show that they agree on $\mathcal{E}'(\Heis)$.
By Lemma \ref{lem:sp(V) action on V} the two derivations agree on $\heis$.
Then Lemma \ref{lem:vanishing D} implies they agree on $\mathcal{E}'(\Heis)$.
\end{proof}

The Lie algebra $T^{-1}\cdot S^2V$ is contained in $\Alg_H(\Heis)$. Therefore we have
\[ \nu:\mathfrak{sp}(V)\to \Alg_H(\Heis)\]
By the way convolution of homogeneous distributions in $\cR_q(\Heis)$ is defined,
Proposition \ref{prop:sp(V) action is inner}
also applies to $\Alg_H(\Heis)$.

\begin{corollary}\label{cor:sp(V) action inner}
The action of $\mathfrak{sp}(V)$ on $\Alg_H(\Heis)$ is inner,
\[ A.f=[\nu(A),f]\qquad A\in \mathfrak{sp}(V),\;f\in\Alg_H(\Heis)\]
\end{corollary}

\subsection{Curvature}\label{sec:curvature}

For a closed contact manifold $M$,
$\cR^j(M)$ denotes the space of smooth sections in the fiber bundle with fiber
the Fr\'echet space $\cR^j(\Heis_x)$.
A smooth section is defined by a local trivialization of the bundle of Heisenberg groups $T_HM$.
$\Alg_H(M)$ is the union
\[ \Alg_H(M) = \bigcup_{k} \sum_{j=0}^k \Tsf^{-j} \cR^{2j}(M)\]
Curvature is an inner derivation of $\Alg_H(M)$.
We briefly recall the construction and refer to \cite{GvE20} for details.

Fix a unitary connection $\nabla$ on $H^{1,0}$,
represented in a local unitary frame as $\nabla=d+\beta$ with a connection $1$-form
\[ \beta\in \Omega^1(M,\mathfrak u(H^{1,0}))\]
The curvature of $\nabla$ is a $\mathfrak u(H^{1,0})$-valued $2$-form
\[ \theta:= d\beta+\frac{1}{2}[\beta,\beta]\in \Omega^2(M, \mathfrak u(H^{1,0}))\]
Let $P_H$ be the principal $U(n)$-bundle of orthonormal frames of $H^{1,0}$.
The unitary group $U(n)$ acts by automorphisms on the model algebra $\Alg_H=\Alg_H(\Heis)$ for $\Heis=\CC^n\times\RR$.
Identify $\Alg_H(M)$ with the algebra of smooth sections of the associated bundle,
\[P_H\times_{U(n)}\Alg_H\]
Viewing $\Alg_H(M)$ as an associated bundle, $\nabla$ determines a covariant derivative
\[ \nabla:\Alg_H(M)\to \Omega^1(\Alg_H(M))\qquad \nabla(a):=da+\beta\cdot a\]
This covariant derivative extends to a derivation on the algebra of $\Alg_H(\Heis_x)$-valued differential forms,
\begin{equation*}
\con  \colon  \Omega^k(\Alg_H(M)) \to \Omega^{k+1}(\Alg_H(M))
\end{equation*}
Applying the map $\nu$ in each fiber
\[ \nu:\mathfrak{sp}(H^{1,0}_x)\to \Alg_H(\Heis_x)\]
we obtain an $\Alg_H(\Heis_x)$-valued 2-form representing curvature,
\[ \curv := \nu(\theta)\in \Omega^2(\Alg_H(M))\]
with
\begin{equation*}
    \con^2(\eta) = [\curv,\eta]\qquad  \eta\in \Omega^\bullet(\Alg_H(M))
\end{equation*}
We have the Bianchi identity
\[ \con(\curv)=0\]
The triple $(\Alg_H(M), \con, \curv)$ is a curved differential graded algebra.

\section{Orbital traces}
\label{sec:orbital-traces}

\subsection{Twisted orbital traces}
Assume given a unitary $\gamma \in U(V)$.
$\gamma$ defines an automorphism of $\Heis$, also denoted $\gamma$, by
\[\gamma(z, t):= (\gamma z, t)\qquad (z,t)\in V\times \RR = \Heis\]
It determines an automorphism of the convolution algebra $C_c^\infty(\Heis)$,
\[f^\gamma(g):=f(\gamma^{-1}g)\qquad f\in C_c^\infty(\Heis)\]
A linear map $\psi:C_c^\infty(\Heis)\to \CC$
is a \emph{$\gamma$-twisted trace} if
\[ \psi(f*g)=\psi(g*f^\gamma) \]
$\gamma$-twisted conjugation is a smooth action of $\Heis$ (as a group) on $\Heis$ (as a manifold),
\[ \mathrm{Ad}_\gamma(a): \Heis\to \Heis\qquad h \mapsto a h \gamma(a^{-1})\]
The $\gamma$-twisted conjugacy class of $h\in \Heis$ is an embedded submanifold of $\Heis$,
\[ [h]_\gamma = \{a h \gamma(a^{-1})\mid a\in \Heis\}\]
The stabilizer of $h\in \Heis$ under the $\mathrm{Ad}_\gamma$ action is the $\gamma$-twisted centralizer,
\[
    Z_\gamma(h) =\{a \in \Heis\mid a h = h \gamma(a) \}
\]
The orbit $[h]_\gamma$ naturally identifies with the homogeneous space $\Heis/Z_\gamma(h)$.
Since the Heisenberg group is nilpotent, both $\Heis$ and $Z_\gamma(h)$ are unimodular.
Hence, by Weil's criterion, $\Heis/Z_\gamma(h)$ admits an $\Heis$-invariant Radon measure, unique up to scaling.
With such an invariant measure, the orbital integral
\[ \psi(f) = \int_{[h]_\gamma} f(z)\,d\mu(z)\qquad f\in C_c^\infty(\Heis)\]
defines a $\gamma$-twisted trace for the convolution algebra $C_c^\infty(\Heis)$.

\subsection{Orbital traces for central elements}\label{sec:twisted orbital traces on Heis}
As above, let $\gamma \in U(V)$.
Decompose $V$ into fixed and moving subspaces,
\[
      V=F\oplus N      \qquad F=\ker(1-\gamma)\qquad N=F^\perp
\]
and
\[
      r=\dim_{\CC}F\qquad m=\dim_{\CC}N\qquad n=r+m
\]
The restriction $\gamma_N=\gamma|_N$ has no fixed vectors, and $\gamma_N-1$ is invertible.

\begin{proposition}
For $c\in \RR$, the $\gamma$-twisted conjugacy class of $(0,c)\in \Heis$ is
\[ \Gamma_{c} = \{(z,c+q(z))\in\Heis \mid z\in N\}\]
where $q:N\to \RR$ is the quadratic form on $N$ defined by
\[q(z)=\ang{Cz,z} \qquad C=\frac{i}{4}(\gamma_N+1)(\gamma_N-1)^{-1}\]
Under the displayed parametrization, Lebesgue measure on $N$ is invariant under $\gamma$-twisted conjugation.
\end{proposition}
\begin{proof}
For twisted conjugation by $(w,s)$, the central coordinate $s$ has no effect and
\begin{equation*}
(w,0)\cdot (0,c)\cdot\gamma (-w,0)
    = ((1-\gamma) w, c+\frac{1}{2}\omega(w,\gamma w))
\end{equation*}
Thus, $(w,s)$ acts trivially iff $w\in F$, and the twisted centralizer is
$Z_\gamma((0,c))=F\times \RR$.
The twisted conjugacy class $\Gamma_c$ identifies with $\Heis/Z_\gamma = N$ via
\[ N\cong \Gamma_c \qquad w\mapsto ((1-\gamma)w, c+\frac{1}{2}\omega(w,\gamma w))\]
With $z=(1-\gamma_N)w$, a direct calculation gives (for $w\in N$)
\[
\frac{1}{2}\omega(w,\gamma w)=\ang{Cz,z}
\]
So we have
\[ N\cong \Gamma_c \qquad z\mapsto (z, c+q(z))\]
Twisted conjugation acts by translation on the $N$-coordinate and therefore preserves Lebesgue measure.
\end{proof}

\begin{definition}
The $\gamma$-twisted trace $\orbtr{c}:C_c^\infty(\Heis)\to \CC$ is integration over the $\gamma$-twisted conjugacy class $\Gamma_c$,
\[\orbtr{c}(f) := \int_N f(z,c+q(z))\,dz\]
Lebesgue measure on $N$ is normalized by $dz=|\omega^m|/m!$.
\end{definition}
For two functions $f,g\in C_c^\infty(\Heis)$ we have
\[ \orbtr{c}(f\ast g) = \orbtr{c}(g\ast f^\gamma)\]
where $f\ast g$ is convolution product, and $f^\gamma(w,s)=f(\gamma^{-1}w,s)$.

\begin{remark}
In the degenerate case where $\gamma$ acts trivially we have $N=\{0\}$. Then $\Gamma_c=\{(0,c)\}$ consists of a single element and
\[ \orbtr{c}(f)=f(0,c)\]
This is the trace considered in \cite{GvE25}.
\end{remark}

\subsection{Orbital traces of symbols}

We now prove the twisted-orbit analogues of results in \cite[Section 2]{GvE25}.
The proofs of \cite{GvE25} generalize easily.

\begin{lemma}
\label{prop:mixed-twisted-orbit}
Let $u\in\mathcal E'(\Heis)$ and $f\in C_c^\infty(\Heis)$. Then
\begin{equation*}
       \orbtr{c}(u*f)=\orbtr{c}(f*u^\gamma)\qquad \orbtr{c}(f\ast u)=\orbtr{c}(u*f^\gamma)
\end{equation*}
\end{lemma}

\begin{proof}
Fix a compact set $C\subset \Heis$ containing the support of $u$ and a slightly larger compact set $C_1$.
Convolution with the fixed test function $f$ gives two continuous maps
\[
 \mathcal E'_{C_1}(\Heis)\to C_c^\infty(\Heis)
 \qquad
 w\mapsto w*f
 \qquad
 w\mapsto f*w^\gamma
\]
whose images have supports in fixed compact sets. Restriction of these
smooth functions to $\Gamma_{c}$, followed by integration there, is
continuous: only the intersection of the orbit with those fixed compact
supports is involved. Choose a compactly supported approximate identity
and regularize $u$ to functions
$u_\varepsilon\in C_c^\infty(\Heis)$ supported in $C_1$, with
$u_\varepsilon\to u$ in $\mathcal E'_{C_1}(\Heis)$.
Then
\[ \orbtr{c}(u*f)=\lim_{\varepsilon\downarrow 0} \orbtr{c}(u_\varepsilon\ast f)=\lim_{\varepsilon\downarrow 0} \orbtr{c}(f\ast u_\varepsilon^\gamma) =\orbtr{c}(f*u^\gamma)\]
The same argument gives the other equality.
\end{proof}

\begin{lemma}
\label{thm:smooth-twisted-trace}
Let $u,v\in\mathcal E'(\Heis)$ be compactly supported regular distributions, i.e.\ smooth away from the identity. Then for $c\ne 0$,
\begin{equation}\label{eq:compact-regular-twisted-trace}
       \orbtr{c}(u*v)=\orbtr{c}(v*u^\gamma).
\end{equation}
\end{lemma}

\begin{proof}
This is the cutoff argument used in \cite[Lemma~2.3]{GvE25}.
Since the closed orbit $\Gamma_{c}$ misses $(0,0)$,
choose a $\gamma$-invariant neighborhood $\mathcal O$ of $(0,0)$ such that $\mathcal O^2$ does not intersect $\Gamma_{c}$.
Choose a $\gamma$-invariant cut-off function $\phi\in C_c^\infty(\mathcal O)$ which equals one near $(0,0)$. Expand
\[
 u*v=(\phi u)*(\phi v)
 +(1-\phi)u*(\phi v)+(\phi u)*(1-\phi)v
 +(1-\phi)u*(1-\phi)v
\]
The first term vanishes on the orbit because it is supported in
$\mathcal O^2$; so does $(\phi v)*(\phi u)^\gamma$.
In each of the other three terms at least one factor
is smooth and compactly supported. Lemma
\ref{prop:mixed-twisted-orbit} applies to each ordered pair of factors.
Since $\phi^\gamma=\phi$, the three resulting terms are exactly the remaining
three terms in the expansion of $v*u^\gamma$. This proves the lemma.
\end{proof}

\begin{proposition}\label{prop:twisted orbital traces on symbols}
If $c\ne 0$ then
\begin{equation*}
\orbtr{c}:\Symb(\Heis)\to \CC \quad \orbtr{c}\left(u\right) := \int_N u(z, c+q(z))\,dz
\end{equation*}
is a $\gamma$-twisted trace on $\Symb(\Heis)$.
\end{proposition}
\begin{proof}
Note that if $u$ is smooth on $\Heis\setminus \{e\}$ and homogeneous of degree $-Q$,
then the restriction of $u$ to $\Gamma_c$ is smooth and asymptotic to $z^{-Q}$.
Therefore $u|_{\Gamma_c}$ is integrable, and $\orbtr{c}(u)$ is well-defined.
We have
\[ \orbtr{c}(u) = a^{Q-2m}\orbtr{a^2c}(u)\qquad a>0\]
By definition of the product in $\Symb(\Heis)$,
\[ \orbtr{c}(u * v) = \lim_{a\downarrow 0} a^{Q-2m} \orbtr{a^2c}(\phi u * \phi v) \]
which, by Lemma \ref{thm:smooth-twisted-trace}, is equal to
\[ \orbtr{c}(v * u^\gamma) = \lim_{a\downarrow 0} a^{Q-2m} \orbtr{a^2c}(\phi v * \phi u^\gamma) \]
\end{proof}

\subsection{Orbital traces on the enlarged algebra}

The twisted traces $\orbtr{c}$ extend further from $\Symb(\Heis)$ to $\Alg_H(\Heis)$.
The integral
\[\orbint{c}(f) = \int_N f(z, c+q(z))\,dz\]
converges if $f\in \cR^{2j}(\Heis)$.
The restriction of $f$ to the orbit $\Gamma_c$ is smooth
and decays as $|z|^{-Q-2j}$, with $Q+2j=2n+2+2j$.
Since the orbit $\Gamma_c$ has real dimension $\dim_\RR N=2m\le 2n$,
the integral converges.

\begin{definition}\label{def:tau_c enlarged}
For $\Tsf^{-j} f_j\in \Alg_H(\Heis)$ with $f_j\in \cR^{2j}(\Heis)$ let
\[ \orbtr{c}(\Tsf^{-j}f_j):=\frac{r!\,c^j}{(r+j)!} \,\int_N f_j(z, c+q(z))\,dz\]
where $r=\dim_\CC \ker(1-\gamma)$.
\end{definition}

We need to show that Definition \ref{def:tau_c enlarged} is compatible with the defining relation of the balanced tensor product (\ref{eq:balanced tensor product}).

\begin{lemma}\label{lem:tau_c of del_t f}
For $f\in \cR^{2j}(\Heis)$ we have
\[ \orbtr{c}(\Tsf^{-j-1}(-\partial_t f))=\orbtr{c}(\Tsf^{-j}f)\]
\end{lemma}
\begin{proof}
We need to show that
\[ \orbint{c}(-\partial_t f)=\frac{r+j+1}{c}\orbint{c}(f)\]
If $f\in \cR^{2j}(\Heis)$ and $a>0$ we have
\begin{align*}
      \orbint{ac}(f)
    &=\int_N f(z, ac+q(z))\,dz
    =a^{m-Q/2-j}\int_N f(w, c+q(w))\,dw\\
    &=a^{-r-j-1}\orbint{c}(f)
\end{align*}
The scaling law in $c$ gives
\[ \orbint{c}(\partial_t f)=\int_N \frac{\partial}{\partial c}f(z, c+q(z))\,dz=\frac{d}{dc}\orbint{c}(f)=-\frac{r+j+1}{c}\orbint{c}(f)\]
Differentiation under the integral is justified, since for a distribution $f\in\cR^{2j}(\Heis)$ the restriction to orbits $\Gamma_c$ is a smooth map
\[ \RR\setminus\{0\}\to L^1(N)\qquad c\mapsto f|_{\Gamma_c}\]
\end{proof}

The extension of $\orbtr{c}$ to $\Alg_H(\Heis)$ has the further property that
\[ \orbtr{c}(\Tsf^{-j}f)=0\qquad \text{for all } f\in \mcU^{2j}(\heis)\]
Note that $\mcU(\heis)$ is precisely the set of distributions supported at the identity.

The action of  $\gamma$ on $\Symb(\Heis)$ fixes the central variable.
The action of $\gamma$ on $\Alg_H(\Heis)$ is defined by
\[ (\Tsf^{-j}f)^\gamma:=\Tsf^{-j}f^\gamma\]
Then $\orbtr{c}$ (with $c\ne 0$)  is a $\gamma$-twisted trace on $\Alg_H(\Heis)$.

\subsection{Topology of the symbol algebra}\label{sec:topology}

A homogeneous distribution $f\in \cR^j(\Heis)$
is determined by its restriction to the sphere $\Sigma_H$,
\[ \Sigma_H = \{(z,t)\in \Heis\mid \|z\|^4+t^2=1\}\]
modulo a distribution supported at the identity $e\in \Heis$.
The latter identifies with an element in $\mcU(\heis)$, and we get
\[ 0\to \mcU^j(\heis)\to \cR^j(\Heis)\to C^\infty(\Sigma_H)\]
Not every smooth function $f\in C^\infty(\Sigma_H)$ extends to a homogeneous distribution of degree $-Q-j$.
Let $\mathcal{P}^j(\Heis)$ be the finite dimensional space of homogeneous polynomials on $\Heis$ of degree $j$ (where $z$ has weight 1 and $t$ has weight 2).
Then $f$ extends if and only if it is orthogonal to $\mathcal{P}^j(\Heis)$,
\[ \int_{\Sigma_H} f(\sigma)p(\sigma)\,d\sigma=0\qquad \forall p\in \mathcal P^j(\Heis)\]
Here $d\sigma$ is the measure on $\Sigma_H$ determined by parabolic polar coordinates,
\[
 \int_{\Heis}\varphi(h)\,dh
 =\int_0^\infty\int_{\Sigma_H}
    \varphi(\delta_r\sigma)r^{Q-1}\,d\sigma\,dr.
\]
For $j=0$ this is the familiar condition that $f$ has mean zero, which is needed for the principal value distribution.
For any function with $f\perp \mathcal{P}^j(\Heis)$ there is a natural choice of homogeneous distribution of degree $-Q-j$ that extends it.\footnote{See \cite[Section 3.2]{HormI} for the analogous situation on $\RR^n$.}

We obtain a decomposition
\[ \cR^{j}(\Heis)=\mathcal{P}^j(\Heis)^\perp\oplus \mcU^j(\heis)\]
Then $\cR^j(\Heis)$ is a Fr\'echet space,
where the finite codimension subspace $\mathcal{P}^j(\Heis)^\perp\subset C^\infty(\Sigma_H)$ has the standard $C^\infty$ topology,
while $\mcU^j(\heis)$ is finite dimensional.

Note that if $j<k$ then $\Tsf^{-j}f=\Tsf^{-k}(\Tsf^{k-j}*f)$, so
\[\Tsf^{-j}\cR^{2j}(\Heis)\subset \Tsf^{-k}\cR^{2k}(\Heis)\]
The algebra $\Alg_H(\Heis)$ is an increasing  union of Fr\'echet spaces,
\[ \Alg_H(\Heis) = \bigcup_{j} \Tsf^{-j} \cR^{2j}(\Heis)\]
Convolution is continuous (see \cite{CGGP92}), 
\[ \cR^{2j}(\Heis)\times\cR^{2k}(\Heis)\to\cR^{2j+2k}(\Heis)\]

\subsection{Continuity of the trace}

Recall that a homogeneous distribution $f\in \cR^{2j}(\Heis)$
is determined by its restriction to the sphere $\Sigma_H$
modulo a distribution supported at $e\in \Heis$.

\begin{lemma}\label{lem:tau_c continuous}
$\orbtr{c}$ is continuous as a map $\Tsf^{-j}\cdot \cR^{2j}(\Heis) \to \CC$ for the sup norm on $C^\infty(\Sigma_H)$,
\[ \orbtr{c}(\Tsf^{-j}f)\le C \|f|_{\Sigma_H}\|_\infty \qquad f\in \cR^{2j}(\Heis)  \]
\end{lemma}
\begin{proof}
Denote by $\rho(z,t)=(|z|^4+t^2)^{1/4}$ the homogeneous norm on $\Heis$,
so that $\Sigma_H=\{h\mid \rho(h)=1\}$.
Then
\[ |\tau_c(\Tsf^{-j}f)| \le \frac{r!|c|^j}{(r+j)!}\, \|f|_{\Sigma_H}\|_\infty \int_N \rho(z,c+q(z))^{-Q-2j}\,dz\]
The integral is finite. 

\end{proof}

\subsection{Contact manifolds}\label{sec:orbital traces on M}

The action of the compact Lie group $G$ on $M$ preserves the contact form $\contact$.
The complex structure $J$ of $H=\ker \contact$ is $G$-invariant.
Fix $\gamma\in G$.
At a point $x\in M^\gamma$ in the fixed-point set
the action of $\gamma$ on $H^{1,0}_x$ is unitary,
\[ \gamma\in U(H^{1,0}_x)\]
Using this unitary, the orbital traces of section \ref{sec:orbital-traces}
are defined at each $x\in M^\gamma$,
\[\orbtr{c}:\Alg_H(\Heis_x)\to \CC\]
Choose a $G$-invariant unitary connection for $\Hzo$, i.e.
\[\gamma(\nabla_X s)=\nabla_{\gamma X}(\gamma s)\qquad (\gamma s)(y)=\gamma( s(\gamma^{-1}y))\]
for sections $s$ of $\Hzo$.
For vector fields $X$ that are tangent to $M^\gamma$ we have $\nabla_X(\gamma s)=\gamma\nabla_X s$,
which means that the endomorphism $\gamma$ of $\Hzo$ is parallel along $M^\gamma$.
Parallel transport along paths in $M^\gamma$ therefore conjugates the values of $\gamma$ and maps eigenspaces to eigenspaces.
Consequently the eigenvalues of $\gamma$ on $\Hzo_x$, with their multiplicities, are constant on each connected component of $M^\gamma$, and the (mutually orthogonal) eigenbundles are smooth.
Thus, the restriction of $\Hzo$ to $M^\gamma$ has local unitary trivializations in which $\gamma$ is represented by a constant diagonal matrix.

With such a trivialization, an element in $\Alg_H(M^\gamma)$ is represented, locally, as a smooth function
\[ M^\gamma\supset U\to \Alg_H(\Heis)\]
for an open set $U$ and
fixed $\Heis=\CC^n\times \RR$.
Here $\Alg_H(M^\gamma)$ denotes sections of the ambient symbol bundle restricted to $M^\gamma$. Its fibers include the normal directions.

In general, $\orbtr{c}$ depends on the action of $\gamma$ on $H^{1,0}_x$, and so a priori $\orbtr{c}$ varies with $x\in M^\gamma$.
But the action of $\gamma$ is constant in our trivialization. Thus $\orbtr{c}$ is the same in each trivialized fiber,
\[ \orbtr{c}:\Alg_H(\Heis)\to \CC\]
Continuity of $\orbtr{c}$ then implies that the composition
\[ U\to \Alg_H(\Heis)\stackrel{\orbtr{c}}{\longrightarrow} \CC\]
is a smooth function.
Thus we get a twisted trace
\[ \orbtr{c}:\Alg_H(M^\gamma)\to C^\infty(M^\gamma)\]
The twisted trace $\orbtr{c}$ is a trace when restricted to $\gamma$-invariant elements.

\begin{definition}
For $c\ne 0$, the trace
\[ \orbtr{c}:\Alg_H(M)^G\to C^\infty(M^\gamma) \]
on the $G$-invariant part of the enlarged symbol algebra $\Alg_H(M)$ is
obtained by restriction to $M^\gamma$ composed with
\[\orbtr{c}:\Alg_H(\Heis_x)\to \CC\]
at each $x\in M^\gamma$.
\end{definition}

\section{The character map}\label{sec:character map}

In this section  we use the closed graded traces $\orbtr{c}$ for the curved dga $(\Alg_H(M), \con,\curv)$ to define the character map
\[ \chi_\gamma:K_1^G(\Symb)\to H^{\mathrm{odd}}(M^\gamma;\CC)\]

\subsection{The traces $\Tr_\gamma$}\label{sec:Tr_gamma}

In \cite{GvE25} we showed that the trace $\tau$ used in the character formula (\ref{eq:chi non-equivariant}) is, in each fiber $\Symb(\Heis_x)$, a linear combination
\footnote{The formula in \cite{GvE25} differs from the one given here by a sign. The difference is a result of changed sign conventions.}
\[ \tau = \frac{(2\pi i)^{n+1}}{2(n!)}\Bigl(\orbtr{+1}-(-1)^n\orbtr{-1}\Bigr)\]
Here $\orbtr{+1}$ and $\orbtr{-1}$ are the traces obtained by evaluation at central elements $(0,c)\in \Heis_x = H_x\times\RR$ with $c=\pm 1$.
We now generalize this to the equivariant case.

Recall that $G$ preserves $\omega$ and the complex structure $J$ (Section~\ref{sec:contact}).
Hence $G$ preserves the Hermitian structure of $\Hzo$.
Choose a $G$-invariant unitary connection $\nabla$ on $\Hzo$.

\begin{definition}\label{def:Tr_gamma}
The trace
\[ \Tr_\gamma:\Alg_H(M)^G\to C^\infty(M^\gamma)\]
is
\[ \Tr_\gamma = \kappa_\gamma \cdot \frac{(2\pi i)^{r+1}}{2(r!)}\Bigl(\orbtr{+1}-(-1)^r\orbtr{-1}\Bigr)\]
where $\orbtr{+1}, \orbtr{-1}$ are the traces $\orbtr{c}$ with $c=\pm 1$ (Section~\ref{sec:orbital traces on M}),
and $r$ and $\kappa_\gamma$ are constants on each connected component of $M^\gamma$,
\[ r=\dim_\CC \ker (1-\gamma)\qquad \kappa_\gamma=\prod_{j=1}^{n-r} \frac{i}{2 \sin (\vartheta_j/2)}\qquad \kappa_\gamma=1\text{ if } r=n\]
Here $e^{i\vartheta_j}$ are the eigenvalues of $\gamma$ (other than $1$)
with $0<\vartheta_j<2\pi$.
\end{definition}

\begin{remark}
Let $\gamma_N^{1/2}$ be the square root of $\gamma_N$ with eigenvalues $e^{i\vartheta_j/2}$, i.e.\ all eigenvalues of $\gamma_N^{1/2}$ have positive imaginary part. Then
\[ \kappa_\gamma = \det\left(\frac{\gamma_N^{1/2}}{1-\gamma_N}\right)\]
If $\gamma$ acts trivially then $M^\gamma=M$, $r=n$, $\kappa_\gamma=1$ and $\Tr_\gamma$ agrees with $\tau$ on $\Symb(M)$.
\end{remark}

The functionals $\orbtr{\pm1}$ are $\gamma$-twisted traces on $\Alg_H(\Heis_x)$, and the restriction to $M^\gamma$ of a $G$-invariant element of $\Alg_H(M)$ is $\gamma$-invariant in each fiber.
Therefore $\Tr_\gamma$ is a trace on $\Alg_H(M)^G$.

\subsection{Differential forms and the connection}\label{sec:Tr_gamma forms}

Let $\con$ and $\curv$ be defined by $\nabla$ as in Section~\ref{sec:curvature}.
Then $\con$ commutes with the action of $G$ on $\Omega^\bullet(\Alg_H(M))$.
The curvature $\theta$ of $\nabla$ is $G$-invariant, and so is $\curv=\nu(\theta)$.
Thus $\con$ is a graded derivation of the algebra  of $G$-invariant forms, 
\[ \nabla: \Omega^\bullet(\Alg_H(M))^G\to \Omega^\bullet(\Alg_H(M))^G\]
and
\[ \con^2(\eta)=[\curv,\eta]\qquad \con(\curv)=0\qquad \curv\in\Omega^2(\Alg_H(M))^G\]
The traces extend to $G$-invariant forms in the usual way.
For  $\eta\in\Omega^k(\Alg_H(M))^G$ and tangent vectors $X_1,\dots,X_k\in T_xM^\gamma$, we set
\[\orbtr{c}:\Omega^k(\Alg_H(M))^G\to\Omega^k(M^\gamma)\qquad \orbtr{c}(\eta)(X_1,\dots,X_k):=\orbtr{c}\bigl(\eta(X_1,\dots,X_k)\bigr)\]
Then $\tau_c$ is a graded trace.

\begin{lemma}\label{lem:Tr_gamma closed}
The trace $\tau_c$ (with $c\ne 0$) is closed,
\[ \orbtr{c}(\con\eta)=d\,\orbtr{c}(\eta)\qquad \eta\in\Omega^k(\Alg_H(M))^G\]
\end{lemma}
\begin{proof}
For $\gamma=e$ this is \cite[Section~5.7]{GvE20}. Locally the connection is
\[ \con=d+[\nu(\beta),\cdot\,]\]
and the trace vanishes on commutators.
In the equivariant case we check that the commutator term is a commutator of $\gamma$-invariant elements.

The statement is local on $M^\gamma$.
Choose $U\subset M^\gamma$ and a unitary frame of $\Hzo|_U$  so that $\gamma$ is a constant diagonal matrix, parallel for the connection,
\[ [\beta,\gamma]=0\]
In this frame, the restriction of $\eta$ to $U$ is a form on $U$ with values in the fixed algebra $\Alg_H(\Heis)$, $\Heis=\CC^n\times\RR$, and on $U$ we have
\[\con\eta=d\eta+[\nu(\beta),\eta]\]
with the graded commutator.
Since $\beta$ commutes with $\gamma$ and $\nu$ is equivariant, $\nu(\beta)$ is a $\gamma$-invariant $1$-form with values in $\Alg_H(\Heis)$,
\[ \nu(\beta)\in \Omega^1(\Alg_H(\Heis))^\gamma\]
The coefficients of $\eta$ are $\gamma$-invariant as well.
Since  $\orbtr{c}$ vanishes on graded commutators of $\gamma$-invariant elements we get 
\[\orbtr{c}([\nu(\beta),\eta])=0\]
It remains to show $\orbtr{c}(d\eta)=d\orbtr{c}(\eta)$.
In our frame, $\orbtr{c}$ is one and the same linear functional on $\Alg_H(\Heis)$ at every point of $U$.
The coefficients of $\eta$ lie in  $\Tsf^{-j}\cR^{2j}(\Heis)$ for fixed $j$ and are smooth on $U$.
On this Fr\'echet space $\orbtr{c}$ is continuous, and so it commutes with differentiation in the directions of $U$.

\end{proof}

\subsection{The character $\chi_\gamma$}\label{sec:chi_gamma}

The character $\chi$ of \cite{GvE20} is the character of a generalized cycle in the sense of \cite{gor1,gor2}.
We refer to \cite[Sections~6.2 and 6.3]{GvE20} for the construction.

Here we have the following.
\begin{enumerate}
\item[(a)] The curved dga
\[(\Omega^\bullet(\Alg_H(M))^G,\con,\curv)\]
where  $\con$ is a graded derivation of degree $1$, $\con^2=[\curv,\cdot\,]$ and $\con(\curv)=0$.
\item[(b)] The closed graded trace
\[\Tr_\gamma:\Omega^\bullet(\Alg_H(M))^G\to\Omega^\bullet(M^\gamma)\qquad \]
into the dga $(\Omega^\bullet(M^\gamma),d)$, with $\Tr_\gamma\circ\con=d\circ\Tr_\gamma$.
\item[(c)] The inclusions
\[ \Symb(M)^G\subset\Alg_H(M)^G=\Omega^0(\Alg_H(M))^G\]
which are  unital homomorphisms.
\end{enumerate}
The morphism of complexes of \cite[Theorem~6.1]{GvE20}  is formulated for an algebra over $C^\infty(M)$ with a $C^\infty(M)$-linear trace into $C^\infty(M)$.
Proofs use only properties (a), (b) and (c), and so apply here, with $\Omega^\bullet(M^\gamma)$ in place of $\Omega^\bullet(M)$.

Moreover, if $\boldsymbol\rho\in\Omega^1(\Alg_H(M))^G$, then $\con'=\con+[\boldsymbol\rho,\cdot\,]$ and $\curv'=\curv+\con(\boldsymbol\rho)+\boldsymbol\rho^2$ again satisfy (a) and (b) with the same trace $\Tr_\gamma$, since $\Tr_\gamma$ vanishes on graded commutators of invariant forms.
The morphisms defined by $(\con,\curv)$ and $(\con',\curv')$ are chain homotopic \cite[Proposition~6.3]{GvE20}.

\begin{definition}\label{def:chi_gamma}
Let $\sigma\in\Symb(M)^G$ be invertible.
The odd form $\chi_\gamma(\sigma)\in\Omega^{\mathrm{odd}}(M^\gamma)$ is
\begin{multline*}
\chi_\gamma(\sigma) = \sum \limits_{l\ge0} \sum \limits_{i_0, \ldots, i_{2l+1} \ge 0} \left(\frac{-1}{2 \pi i}\right)^{I+l+1} \frac{l!}{(I+2l+1)!}\\
 \Tr_\gamma \left(  \sigma^{-1} \curv^{i_0} \con(\sigma)\curv^{i_1}\con(\sigma^{-1}) \dots \con(\sigma) \curv^{i_{2l+1}}\right)
\end{multline*}
\[I=i_0+i_1+\cdots+i_{2l+1}\]
\end{definition}
The sum is finite: The term with indices $l$, $i_0,\dots,i_{2l+1}$ has form degree $2I+2l+1$,
and on a component of $M^\gamma$ of dimension $2r+1$ only terms with $I+l\le r$ contribute.

The formula is that  of \cite[Theorem~6.7]{GvE20}, with $\tau$ replaced by $\Tr_\gamma$.
For $\gamma=e$, the trace $\Tr_e$ agrees with $\tau$ on $\Symb(M)$, and $\chi_e$ is given by the same formula as the character $\chi$ of \cite{GvE20}, now in the kernel picture.

\subsection{The character on equivariant $K$-theory}\label{sec:chi_gamma well defined}

Since $\Symb(M)$ is closed under holomorphic functional calculus in its $C^*$-algebra completion \cite[Appendix]{GvE20},
$K_1^G(\Symb)$ is  unchanged if $\Symb(M)$ is replaced by its $C^*$-algebra completion.

Let $W$ be a finite-dimensional unitary representations of $G$.
The group $G$ acts on $\Omega^\bullet(\Alg_H(M))\otimes\End W$, on $\End W$ by conjugation.
We extend $\con, \curv$ by $\con\otimes1, \curv\otimes1$, and $\orbtr{c}$ by the matrix trace twisted by $\gamma$, 
\[\orbtr{c}(a\otimes m):=\orbtr{c}(a)\,\tr_W(m\gamma)\]
This gives $\Tr_\gamma$ on $(\Omega^\bullet(\Alg_H(M))\otimes\End W)^G$.
Sections~\ref{sec:Tr_gamma forms} and \ref{sec:chi_gamma} go through verbatim with these coefficients, and  (a)--(c) hold. 
Definition~\ref{def:chi_gamma} of $\chi_\gamma(\sigma)$ applies to an invertible $\sigma\in(\Symb(M)\otimes\End W)^G$.


\begin{theorem}\label{thm:chi_gamma well defined}
Let $W$ be a finite-dimensional unitary representation of $G$, and let $\sigma\in(\Symb(M)\otimes\End W)^G$ be invertible.
The form $\chi_\gamma(\sigma)$ is closed.
Its cohomology class depends only on the $K$-theory class  $[\sigma]\in K_1^G(\Symb)$, and is independent of the choice of $\nabla$ or $J$.
We obtain a homomorphism
\[ \chi_\gamma:K_1^G(\Symb)\to H^{\mathrm{odd}}(M^\gamma;\CC)\]
\end{theorem}
\begin{proof}
Once everything is $G$-invariant, the proof is that of \cite[Theorem~6.7]{GvE20}.

\emph{Closedness and homotopy invariance.}
By Section~\ref{sec:chi_gamma}, \cite[Theorem~6.1]{GvE20} gives a morphism of complexes from the periodic cyclic complex of $(\Symb(M)\otimes\End W)^G$ to $(\Omega^\bullet(M^\gamma)[u^{-1},u],ud)$.
Since $\Ch(\sigma)$ is a cycle, $\chi_\gamma(\sigma)$ is closed.
The morphism is continuous for the Fr\'echet topology, because $\orbtr{c}$ is continuous (Lemma~\ref{lem:tau_c continuous}), and $\con$ and the product are continuous on finite stages (Section~\ref{sec:topology}).
As in \cite[Section~6.3]{GvE20}, the cohomology class of $\chi_\gamma(\sigma)$ therefore depends only on the homotopy class of $\sigma$ among invertible elements of $(\Symb(M)\otimes\End W)^G$.

\emph{Relations in $K_1^G(\Symb)$.}
The connection, the curvature and the trace are compatible with direct sums $W\oplus W'$, so $\chi_\gamma(\sigma\oplus\sigma')=\chi_\gamma(\sigma)+\chi_\gamma(\sigma')$.
Since $\con(1_W)=0$ and every term in Definition~\ref{def:chi_gamma} contains a factor $\con(\sigma)$, we have $\chi_\gamma(1_W)=0$.
A $G$-equivariant unitary $U:W\to W'$ is constant, so it commutes with $\con$ and $\curv$, and $\tr_{W'}(UmU^*\gamma)=\tr_W(m\gamma)$ because $U\gamma=\gamma U$.
Hence $\chi_\gamma(U\sigma U^*)=\chi_\gamma(\sigma)$.
This shows that $[\sigma]\mapsto[\chi_\gamma(\sigma)]$ is a well defined homomorphism on $K_1^G(\Symb)$.

\emph{Independence of the connection.}
Let $\nabla_0$, $\nabla_1$ be $G$-invariant unitary connections on $\Hzo$, and $\rho=\nabla_1-\nabla_0\in\Omega^1(M,\Lu(\Hzo))^G$.
With  $\nu(\rho)\in\Omega^1(\Alg_H(M))^G$ as in \cite[Lemma~4.5]{GvE20},
\[ \con_1=\con_0+[\nu(\rho),\cdot\,]\qquad \curv_1=\curv_0+\con_0(\nu(\rho))+\nu(\rho)^2\]
The first identity follows from $\con\eta=d\eta+[\nu(\beta),\eta]$ in a local frame.
The second follows from $\theta_1=\theta_0+\nabla_0\rho+\frac12[\rho,\rho]$, because $\nu$ is a Lie algebra homomorphism.
Locally,
\[\nu(\nabla_0\rho)=\nu(d\rho+[\beta_0,\rho])=d\nu(\rho)+[\nu(\beta_0),\nu(\rho)]=\con_0(\nu(\rho))\] 
and $\nu(\frac12[\rho,\rho])=\nu(\rho)^2$.
By Section~\ref{sec:chi_gamma}, the two characters are chain homotopic, and $[\chi_\gamma(\sigma)]$ does not depend on the connection.

\emph{Independence of $J$.}
Let $J_0$, $J_1$ be $G$-invariant complex structures on $H$ compatible with $\omega$.
They are connected by a path $J_s$ of such structures.
The algebras $\Symb(M)$ and $\Alg_H(M)$, the twisted conjugacy classes $\Gamma_c$ and their measures $|\omega^m|/m!$ depend only on $\omega$ and $\gamma$, not on $J_s$.
The rank $r$ and the eigenvalues $e^{i\vartheta_j}$ of $\gamma$ on $(N,J_s)$ vary continuously with $s$ and lie in a finite set, so they are constant along the path.
Hence $\Tr_\gamma$ is the same for $J_0$ and $J_1$.
Choose $G$-invariant connections $\nabla_0$, $\nabla_1$ that are unitary for $J_0$ and $J_1$ respectively.
Both preserve $\omega$, so $\rho=\nabla_1-\nabla_0\in\Omega^1(M,\spl(H))^G$.
The identities of the previous paragraph hold with $\nu:\spl(H)\to\Alg_H(M)$, since the action of all of $\spl(V)$ is inner (Corollary~\ref{cor:sp(V) action inner}), and the two characters are again chain homotopic.

\end{proof}

\begin{remark}
A $G$-equivariant Hermitian vector bundle $E$ on $M$ is a $G$-invariant summand $E$ of a trivial bundle $M\times W$. Let $p\in C^\infty(M,\End W)^G$ be  the corresponding projection.
If $P$ is a $G$-equivariant Heisenberg elliptic operator acting on sections of $E$, then
\[[\sigma_H(P)]\in K_1^G(\Symb)\]
is, by definition, the class of
\[ \sigma_H(P)+(1-p)\in (\Symb(M)\otimes\End W)^G\]
This is the class to which $\chi_\gamma$ is applied in  Theorem~\ref{thm:main-intro}.
    
\end{remark}

\section{Toeplitz operators}\label{sec:Toeplitz operators}

In this section we compute $\chi_\gamma$ on the symbols of Toeplitz operators.
We first evaluate $\Tr_\gamma$ on the principal symbol $S_+$ of the Szeg\"o projection by an explicit calculation.
The second computation is that $\curv$ acts on the range of $S_+$ as the scalar $2$-form $\frac12\tr(\theta)$.
For Toeplitz operators, the integrand of Theorem~\ref{thm:main-intro} then becomes the integrand of the equivariant index formula of \cite[Theorem~2]{GvE26}.

\subsection{The Szeg\"o symbol}\label{sec:Szego symbol}

The Weyl symbol of the projection onto the ground state of the harmonic oscillator is
\[ s_0(\xi)=2^n e^{-|\xi|^2}\qquad \xi\in V^*\]
Denote Fourier transform on $\Heis$ by
\[\mcF f(\xi,\ft)=\int_{V\times\RR} f(z,t)\,e^{-i\xi(z)-ist}\,dz\,dt\]

\begin{definition}\label{def:Szego symbol}
The \emph{Szeg\"o symbol} $S_+\in\Symb(\Heis)$ is the element whose Weyl pair is
$(s_0,0)$, 
i.e.
\[\mcF S_+(\xi,-1)=s_0(\xi)\qquad \mcF S_+(\xi,1)=0\]
\end{definition}

Since $s_0\#s_0=s_0$ in the Weyl algebra it follows that $S_+$ is a projection in $\Symb(\Heis)$.
It is the principal symbol of the Szeg\"o projection.
Since $s_0$ is invariant under $U(V)$, the elements $S_+\in\Symb(\Heis_x)$, $V=\Hzo_x$, form a $G$-invariant smooth section $S_+\in\Symb(M)$.
A generalized Szeg\"o projection is a projection $S\in\Psi^0_H(M)$ with $\sigma_H(S)=S_+$, and $S$ can be chosen $G$-invariant.

\begin{lemma}\label{lem:Szego kernel}
On $\Heis\setminus\{e\}$ the distribution $S_+$ is the smooth function
\[ S_+(z,t)=\frac{n!}{(2\pi)^{n+1}}\left(\frac{|z|^2}{4}+it\right)^{-n-1}\]
\end{lemma}
\begin{proof}
Since $dz=|\omega^n|/n!$ is Lebesgue measure in unitary coordinates, the inverse Fourier transform on $V$ of $s_0$ is the Gaussian $(2\pi)^{-n}e^{-|z|^2/4}$.
Write
\[\mcF_tf(z,\ft)=\int_\RR f(z,t)\,e^{-i\ft t}\,dt\]
for the partial Fourier transform in the central variable.
Then
\[\mcF_tS_+(z,-1)=(2\pi)^{-n}e^{-|z|^2/4}\]
and homogeneity  gives, for $y>0$,
\[ \mcF_tS_+(z,-y)=\Bigl(\frac{y}{2\pi}\Bigr)^ne^{-y|z|^2/4}\qquad \mcF_tS_+(z,y)=0\]
For $\varepsilon>0$, we insert a factor $e^{-\varepsilon y}$ and invert the Fourier transform in $t$.
The Gamma integral $\int_0^\infty y^ne^{-y\zeta}\,dy=n!\,\zeta^{-n-1}$, $\operatorname{Re}\zeta>0$, gives the bounded smooth function
\[ K_\varepsilon(z,t)=\frac{1}{2\pi}\int_0^\infty\Bigl(\frac{y}{2\pi}\Bigr)^ne^{-y(\varepsilon+|z|^2/4+it)}\,dy=\frac{n!}{(2\pi)^{n+1}}\left(\varepsilon+\frac{|z|^2}{4}+it\right)^{-n-1}\]
By Fubini's theorem, $\ang{K_\varepsilon,\varphi}$ is the integral above with the extra factor $e^{-\varepsilon y}$.
The integrand is dominated by the integrand at $\varepsilon=0$, which is integrable.
Dominated convergence gives $\ang{K_\varepsilon,\varphi}\to\ang{S_+,\varphi}$ as $\varepsilon\downarrow0$.
On the other hand, $K_\varepsilon$ converges uniformly on compact subsets of $\Heis\setminus\{e\}$ to the function in the lemma.
Taking $\varphi$ with compact support in $\Heis\setminus\{e\}$ proves the lemma.

\end{proof}

\subsection{The orbital integrals of $S_+$}\label{sec:Szego orbital integrals}

Let $\gamma\in U(V)$, with $V=F\oplus N$, $r$, $m$, $q$ and $\Gamma_c$ as in Section~\ref{sec:twisted orbital traces on Heis}.
Let $\gamma_N^{1/2}$ be the square root of $\gamma_N$ with eigenvalues $e^{i\vartheta_j/2}$, $0<\vartheta_j<2\pi$, as in Section~\ref{sec:Tr_gamma}, and set
\[ \sqrt{\det\gamma_N}:=\det\bigl(\gamma_N^{1/2}\bigr)=\prod_{j=1}^m e^{i\vartheta_j/2}\]

\begin{proposition}\label{prop:Szego orbital values}
For $c\ne0$,
\[ \orbtr{c}(S_+)=\frac{r!\,\det(1-\gamma_N)}{(2\pi i c)^{r+1}}\]
Consequently, for $\gamma\in G$,
\[ \Tr_\gamma(S_+)=\sqrt{\det\gamma_N}\]
on $M^\gamma$.
\end{proposition}
\begin{proof}
Choose a unitary basis $e_1,\dots,e_m$ of $N$ with $\gamma e_j=e^{i\vartheta_j}e_j$, $0<\vartheta_j<2\pi$, and write $z=\sum_j z_je_j$.
Then $|z|^2=\sum_j|z_j|^2$, and $dz=|\omega^m|/m!$ is the product of the Lebesgue measures of the coordinates $z_j\in\CC$.
On each eigenline
\[ C e_j=\frac{i}{4}\,\frac{e^{i\vartheta_j}+1}{e^{i\vartheta_j}-1}\,e_j=\frac14\cot\frac{\vartheta_j}{2}\,e_j\]
so
\[ q(z)=\frac14\sum_{j=1}^m\cot\frac{\vartheta_j}{2}\,|z_j|^2\]
By Lemma~\ref{lem:Szego kernel}, the integrand of $\orbtr{c}(S_+)$ is $\frac{n!}{(2\pi)^{n+1}}\Lambda(z)^{-n-1}$, where
\[ \Lambda(z)=\frac{|z|^2}{4}+i\bigl(c+q(z)\bigr)=ic+\sum_{j=1}^mb_j|z_j|^2\qquad b_j=\frac14\Bigl(1+i\cot\frac{\vartheta_j}{2}\Bigr)\]
We have $\operatorname{Re}\Lambda(z)=\frac14|z|^2$ and $\Lambda(0)=ic\ne0$.
Thus $\Lambda$ vanishes nowhere on $N$, $|\Lambda(z)|\ge\frac14|z|^2$, and $|\Lambda|^{-n-1}$ is integrable on $N$ because $2n+2>2m$.

Let $\varepsilon>0$.
Since $\operatorname{Re}\bigl(\varepsilon+\Lambda(z)\bigr)=\varepsilon+\frac14|z|^2>0$, the Gamma integral gives $\bigl(\varepsilon+\Lambda(z)\bigr)^{-n-1}=\frac1{n!}\int_0^\infty y^ne^{-y(\varepsilon+\Lambda(z))}\,dy$.
Fubini's theorem applies, since
\[ \int_0^\infty\!\!\int_Ny^ne^{-y(\varepsilon+|z|^2/4)}\,dz\,dy=\int_0^\infty y^ne^{-y\varepsilon}\Bigl(\frac{4\pi}{y}\Bigr)^m dy<\infty\]
With $\int_\CC e^{-yb|\zeta|^2}\,d\zeta=\pi/(yb)$ for $\operatorname{Re}b>0$, applied to $b=b_j$, $\operatorname{Re}b_j=\frac14>0$, and with $n-m=r$, we obtain
\[ \int_N\bigl(\varepsilon+\Lambda(z)\bigr)^{-n-1}dz=\frac{1}{n!}\int_0^\infty y^ne^{-y(\varepsilon+ic)}\prod_{j=1}^m\frac{\pi}{y b_j}\,dy=\frac{r!\,\pi^m}{n!\,b_1\cdots b_m}\,(\varepsilon+ic)^{-r-1}\]
Since $\operatorname{Re}\Lambda\ge0$, we have $|\varepsilon+\Lambda|\ge|\Lambda|$.
By dominated convergence we may let $\varepsilon\downarrow0$, and
\[ \orbtr{c}(S_+)=\frac{n!}{(2\pi)^{n+1}}\int_N\Lambda(z)^{-n-1}dz=\frac{r!\,\pi^m}{(2\pi)^{n+1}\,b_1\cdots b_m}\,(ic)^{-r-1}\]
Finally, $1+i\cot(\vartheta/2)=ie^{-i\vartheta/2}/\sin(\vartheta/2)$ gives
\[ \frac{1}{b_j}=\frac{4\sin(\vartheta_j/2)}{ie^{-i\vartheta_j/2}}=2\bigl(1-e^{i\vartheta_j}\bigr)\qquad \frac{\pi^m}{b_1\cdots b_m}=(2\pi)^m\det(1-\gamma_N)\]
and so $\orbtr{c}(S_+)=r!\det(1-\gamma_N)/\bigl((2\pi)^{r+1}(ic)^{r+1}\bigr)$, as claimed.
(If $m=0$, the products are empty and the formula is $\orbtr{c}(S_+)=S_+(0,c)$.)

For the second statement, fix $x\in M^\gamma$ and apply the first in $\Symb(\Heis_x)$, $V=\Hzo_x$.
Since $\orbtr{-1}(S_+)=(-1)^{r+1}\orbtr{+1}(S_+)$,
\[ \orbtr{+1}(S_+)-(-1)^r\orbtr{-1}(S_+)=2\,\orbtr{+1}(S_+)=\frac{2\,r!\,\det(1-\gamma_N)}{(2\pi i)^{r+1}}\]
By Definition~\ref{def:Tr_gamma} and the formula $\kappa_\gamma=\det\bigl(\gamma_N^{1/2}(1-\gamma_N)^{-1}\bigr)$ of Section~\ref{sec:Tr_gamma},
\[ \Tr_\gamma(S_+)=\frac12\,\kappa_\gamma\,\frac{(2\pi i)^{r+1}}{r!}\cdot\frac{2\,r!\,\det(1-\gamma_N)}{(2\pi i)^{r+1}}=\kappa_\gamma\det(1-\gamma_N)=\det\bigl(\gamma_N^{1/2}\bigr)\]
\end{proof}

\begin{remark}
For $\gamma=e$ we have $N=\{0\}$, and Proposition~\ref{prop:Szego orbital values} gives $\Tr_e(S_+)=1$.
With   $\Tr_e=\tau$ this agrees with $\tau(S_+)=\Tr(s_0)=1$ as in  \cite{GvE25}.
\end{remark}

\subsection{Toeplitz operators}\label{sec:Toeplitz set-up}

Let $S\in \Psi^0_H(M)$ be a $G$-invariant generalized Szeg\"o projection (cf. \cite{EMxx}), with
\[ \sigma_H(S)=S_+ \qquad S=S^*=S^2 \]
Let $W$ be a finite-dimensional unitary representation of $G$.
Let $f:M\to GL(W)$ be smooth and $G$-equivariant, $f(\gamma x)=\gamma f(x)\gamma^{-1}$,
and let $f$ act on $L^2(M;W)$ by multiplication $(\mathcal M_fu)(x)=f(x)u(x)$.
The Toeplitz operator is 
\[ T_f:=(S\otimes1_W)\mathcal M_f\]
acting on the range of  $S\otimes1_W$.
The operator  $\tilde T_f=T_f\oplus1$ on $L^2(M;W)$ is a $G$-equivariant order zero pseudodifferential operator in the Heisenberg calculus.
It has the same index as $T_f$.
The symbol is
\[ \sigma_H(\tilde T_f)=S_+\otimes f+(1-S_+)\otimes 1_W\]

\subsection{Curvature on the vacuum}\label{sec:vacuum curvature}

The key fact is that $S_+$ is annihilated by the anti-holomorphic directions of $\Heis$.

\begin{lemma}\label{lem:dbar Szego}
For every $v\in V$,
\[ (\mathsf v+i\,\mathsf{Jv})*S_+=0\]
where $\mathsf v = -v \delta_e$ (in the notation of Section~\ref{sec:enlarged-symbol-algebra}).
\end{lemma}
\begin{proof}
Since $\mathsf v$ is supported at $e$, $\mathsf v*S_+$ is the ordinary convolution of distributions: $\mathsf v*((1-\phi)S_+)$ is smooth near $e$, and $\mathsf v*S_+$ is already homogeneous, so it is the principal part of $\phi\mathsf v*\phi S_+$ (Section~\ref{sec:model symbols}).
For a distribution $u$ on $\Heis$, the definition of convolution gives $\ang{\mathsf v*u,\varphi}=\ang{u,v^R\varphi}$, where $v^R\varphi(h)=\frac{d}{d\varepsilon}\big|_{\varepsilon=0}\varphi((\varepsilon v,0)h)$.
By the group law (\ref{eq:heisenberg group law}), $v^R=\partial_v-\frac12\omega(v,z)\partial_t$, a vector field that preserves the Haar measure.
Hence
\[ \mathsf v*u=-\partial_vu+\tfrac12\omega(v,z)\,\partial_tu\]
In terms of the partial Fourier transform $\mcF_t$ of Section~\ref{sec:Szego symbol},
\[ \mcF_t(\mathsf v*u)(z,\ft)=\Bigl(-\partial_v+\tfrac{i\ft}{2}\,\omega(v,z)\Bigr)\mcF_tu(z,\ft)\]
Take $u=S_+$.
By the proof of Lemma~\ref{lem:Szego kernel}, $\mcF_tS_+(z,\ft)=0$ for $\ft>0$, and $\mcF_tS_+(z,-y)=(y/2\pi)^ne^{-y|z|^2/4}$ for $y>0$.
This is a locally integrable function, smooth in $z$ for fixed $\ft$, whose $z$-derivatives are locally integrable, so the right-hand side above may be computed pointwise.
Since $\partial_v|z|^2=2g(v,z)$, for $y>0$ we get
\[ \Bigl(-\partial_v-\tfrac{iy}{2}\,\omega(v,z)\Bigr)e^{-y|z|^2/4}=\tfrac y2\bigl(g(v,z)-i\omega(v,z)\bigr)e^{-y|z|^2/4}=\tfrac y2\ang{v,z}\,e^{-y|z|^2/4}\]
Since the Hermitian form is linear in the first variable, $\ang{Jv,z}=i\ang{v,z}$.
Therefore
\[ \mcF_t\bigl((\mathsf v+i\,\mathsf{Jv})*S_+\bigr)(z,-y)=\tfrac y2\bigl(\ang{v,z}+i\cdot i\ang{v,z}\bigr)\Bigl(\frac{y}{2\pi}\Bigr)^ne^{-y|z|^2/4}=0\]
Since $S_+$ is homogeneous, it is a tempered distribution, and the partial Fourier transform is injective on tempered distributions, which proves the lemma.

\end{proof}

\begin{proposition}\label{prop:vacuum curvature}
For $A\in\Lu(V)\subset \mathfrak{sp}(V)$,
\[ \nu(A)S_+=S_+\nu(A)=\frac12\tr(A)\,S_+\]
in $\Alg_H(\Heis)$, where $\tr(A)$ is the trace of $A$ as a complex linear map of $V$.
\end{proposition}
\begin{proof}
Since $s_0$ is $U(V)$-invariant, so is $S_+$, and $A.S_+=0$.
By Corollary~\ref{cor:sp(V) action inner}, $[\nu(A),S_+]=A.S_+=0$.

Both sides of the identity are linear in $A$.
Choose a unitary basis $w_1,\dots,w_n$ of $V$ with $Aw_k=i\lambda_kw_k$, $\lambda_k\in\RR$.
Then $A=\sum_k\lambda_kA_k$, where $A_kv=i\ang{v,w_k}w_k$ is $J$ on $\CC w_k$ and zero on its orthogonal complement, and $\tr(A_k)=i$.
It suffices to show $\nu(A_k)S_+=\frac i2S_+$.

Fix $k$ and write $w=w_k$.
In $\heis$ we have $[w,Jw]=-\omega(w,Jw)\,T=-g(w,w)\,T=-T$, and the real span of $w,Jw$ is symplectically orthogonal to the other $w_l,Jw_l$.
A direct check then shows that $\mu\bigl(\frac12T^{-1}(w^2+(Jw)^2)\bigr)=A_k$ (Lemma~\ref{lem:sp(V) action on V}), and so
\[ \nu(A_k)=\tfrac12\,\Tsf^{-1}\bigl(\mathsf w*\mathsf w+\mathsf{Jw}*\mathsf{Jw}\bigr)\]
Since $X\mapsto\mathsf X$ preserves brackets, $[\mathsf w,\mathsf{Jw}]=-\Tsf$, and
\[ \mathsf w*\mathsf w+\mathsf{Jw}*\mathsf{Jw}=(\mathsf w-i\,\mathsf{Jw})*(\mathsf w+i\,\mathsf{Jw})-i[\mathsf w,\mathsf{Jw}]=(\mathsf w-i\,\mathsf{Jw})*(\mathsf w+i\,\mathsf{Jw})+i\Tsf\]
By Lemma~\ref{lem:dbar Szego}, $(\mathsf w+i\,\mathsf{Jw})*S_+=0$.
Hence $(\mathsf w*\mathsf w+\mathsf{Jw}*\mathsf{Jw})*S_+=i\,\Tsf*S_+$, and
\[ \nu(A_k)S_+=\tfrac12\Tsf^{-1}\bigl(i\,\Tsf*S_+\bigr)=\tfrac i2\,S_+\]

\end{proof}

\begin{remark}\label{rem:vacuum curvature GvE20}
Proposition~\ref{prop:vacuum curvature} is the kernel-picture version of \cite[Lemma~7.3]{GvE20}, where $Q_ks_0=s_0$ for the harmonic oscillator $Q_k$ in the $k$-th complex direction.
\end{remark}

\begin{corollary}\label{cor:vacuum curvature}
In $\Omega^\bullet(\Alg_H(M))$,
\[ \curv S_+=S_+\curv=\tfrac12\tr(\theta)\,S_+\qquad \curv^jS_+=\bigl(\tfrac12\tr(\theta)\bigr)^jS_+\quad j\ge0\]
Consequently, for $a\in\Omega^\bullet(M;\End W)^G$ and $j\ge0$,
\[ \Tr_\gamma\bigl(\curv^jS_+\otimes a\bigr)=\sqrt{\det\gamma_N}\,\bigl(\tfrac12\tr(\theta)\bigr)^j\,\tr_W(a\gamma)\]
on $M^\gamma$, where $\tr_W(a\gamma)$ is restricted to $M^\gamma$.
\end{corollary}
\begin{proof}
At each $x\in M$, $\curv=\nu(\theta)$ with $\theta\in\Lambda^2T^*_xM\otimes\Lu(\Hzo_x)$, and $\nu$ is applied to the coefficients.
Proposition~\ref{prop:vacuum curvature}, applied coefficientwise, gives $\curv S_+=S_+\curv=\frac12\tr(\theta)S_+$.
Since $\frac12\tr(\theta)$ is a scalar $2$-form, it is central in $\Omega^\bullet(\Alg_H(M))$, and the formula for $\curv^jS_+$ follows by induction on $j$.
The element $\curv^jS_+\otimes a=S_+\otimes b$, $b=(\frac12\tr(\theta))^ja$, is $G$-invariant, since $\theta$ and $a$ are.
By definition  $\orbtr{c}$ is computed pointwise on $M^\gamma$ and is linear over forms on $M^\gamma$, so $\orbtr{c}(S_+\otimes b)=\orbtr{c}(S_+)\,\tr_W(b\gamma)$.
Definition~\ref{def:Tr_gamma} and Proposition~\ref{prop:Szego orbital values} give the trace.

\end{proof}

For $\gamma=e$ the last formula is \cite[Corollary~7.4]{GvE20}.

\subsection{The character of a Toeplitz symbol}\label{sec:Toeplitz character}

For a finite-dimensional unitary representation $W$ of $G$ and a smooth $G$-equivariant map $f:M\to GL(W)$, the odd equivariant Chern character is 
\[ \ch_\gamma(f)=-\sum_{l\ge0}\frac{l!}{(2l+1)!\,(2\pi i)^{l+1}}\,\tr_W\bigl((f^{-1}df)^{2l+1}\gamma\bigr)\in\Omega^{\mathrm{odd}}(M^\gamma)\]

\begin{theorem}\label{thm:Toeplitz character}
With $\tilde T_f$  as in Section~\ref{sec:Toeplitz set-up} we have
\[ \chi_\gamma\bigl(\sigma_H(\tilde T_f)\bigr)=\sqrt{\det\gamma_N}\;\ch_\gamma(f)\wedge\exp\Bigl(\tfrac12\cdot\tfrac{i}{2\pi}\tr(\theta)\Bigr)\]
as differential forms.
In cohomology,
\[ \chi_\gamma\bigl([\sigma_H(\tilde T_f)]\bigr)=\sqrt{\det\gamma_N}\;\ch_\gamma(f)\wedge\exp\Bigl(\tfrac12c_1(\Hzo)\Bigr)\in H^{\mathrm{odd}}(M^\gamma;\CC)\]
\end{theorem}
\begin{proof}
The proof is that of \cite[Theorem~7.5]{GvE20}, with $\Trh\circ\tr$ replaced by $\Tr_\gamma$.
Write $\sigma=\sigma_H(\tilde T_f)$, so that
\[ \sigma=S_+\otimes f+(1-S_+)\otimes 1_W\qquad \sigma^{-1}=S_+\otimes f^{-1}+(1-S_+)\otimes 1_W\]
Since $\con S_+=0$ and  $W$ carries the trivial connection, we have
\[ \con\sigma=S_+\otimes df\qquad \con(\sigma^{-1})=S_+\otimes d(f^{-1})\]
Every term of Definition~\ref{def:chi_gamma} contains a factor $\con\sigma$.
By Corollary~\ref{cor:vacuum curvature} each factor $\curv^{i_j}$ that stands next to a factor $S_+\otimes(\cdots)$ can be replaced by the central form $(\frac12\tr(\theta))^{i_j}$, and $\sigma^{-1}(S_+\otimes df)=S_+\otimes f^{-1}df$.
Since $d(f^{-1})=-f^{-1}df\,f^{-1}$,
\[ \sigma^{-1}\curv^{i_0}\con(\sigma)\curv^{i_1}\con(\sigma^{-1})\cdots\con(\sigma)\curv^{i_{2l+1}}=(-1)^l\bigl(\tfrac12\tr(\theta)\bigr)^I\,S_+\otimes(f^{-1}df)^{2l+1}\]
with $I=i_0+\dots+i_{2l+1}$.
By Corollary~\ref{cor:vacuum curvature}, with $j=0$, its trace is 
\[(-1)^l\sqrt{\det\gamma_N}\,(\frac12\tr(\theta))^I\,\tr_W((f^{-1}df)^{2l+1}\gamma)\]
The number of ways to write $I$ as $i_0+\dots+i_{2l+1}$ is $\frac{(I+2l+1)!}{I!\,(2l+1)!}$, so the coefficients of Definition~\ref{def:chi_gamma} give
\begin{multline*}
\chi_\gamma(\sigma)=\sqrt{\det\gamma_N}\Bigl(\sum_{l\ge0}(-1)^l\Bigl(\frac{-1}{2\pi i}\Bigr)^{l+1}\frac{l!}{(2l+1)!}\tr_W\bigl((f^{-1}df)^{2l+1}\gamma\bigr)\Bigr)\\
\wedge\sum_{I\ge0}\frac{1}{I!}\Bigl(\frac{-1}{2\pi i}\cdot\frac12\tr(\theta)\Bigr)^I
\end{multline*}
The first factor is $\ch_\gamma(f)$, and the second is $\exp(\frac12\cdot\frac{i}{2\pi}\tr(\theta))$.
The form $\frac{i}{2\pi}\tr(\theta)$ restricted to $M^\gamma$ represents $c_1(\Hzo|_{M^\gamma})$.

\end{proof}

\subsection{The fixed-point integrand}\label{sec:Toeplitz integrand}

Along $M^\gamma$ we have $\Hzo|_{M^\gamma}=F\oplus N$, with $F=\ker(1-\gamma)$ the complex contact bundle of $M^\gamma$ and $N=N^\gamma$.
The tangent bundle is $TM^\gamma=F\oplus\RR\Reeb$, and $\Td(M^\gamma):=\Td(F)$.
The element $\gamma$ acts on $N^*$ by the inverse transpose.

\begin{proposition}\label{prop:Toeplitz integrand}
With $\tilde T_f$ as in Theorem~\ref{thm:Toeplitz character},
\[ \chi_\gamma\bigl([\sigma_H(\tilde T_f)]\bigr)\wedge\Ahat(M^\gamma)\wedge\frac{1}{\Delta_\gamma(N^\gamma)}=\ch_\gamma(f)\wedge\Td(M^\gamma)\wedge\frac{1}{\ch_\gamma(\lambda_{-1}N^*)}\]
in $H^\bullet(M^\gamma;\CC)$.
\end{proposition}
\begin{proof}
Apart from $\ch_\gamma(f)$ all classes involved are even, so they commute, and $c_1(\Hzo)=c_1(F)+c_1(N)$ on $M^\gamma$.
For the fixed part, $\Ahat(M^\gamma)=\Ahat(F)$, and 
\[e^{c_1(F)/2}\Ahat(F)=\Td(F)\]
For the moving part, let $x_1,\dots,x_m$ be the Chern roots of $N$, with $x_j$ belonging to the eigenvalue $e^{i\vartheta_j}$, $0<\vartheta_j<2\pi$.
This is the same branch used both for $\sqrt{\det\gamma_N}$ and in $\Delta_\gamma$, so the left-hand side does not depend on the choice of branch.
Then 
\[\sqrt{\det\gamma_N}\,e^{c_1(N)/2}=\prod_je^{i\vartheta_j/2}e^{x_j/2}\]
and
\[ \frac{\sqrt{\det\gamma_N}\,e^{c_1(N)/2}}{\Delta_\gamma(N)}=\prod_{j=1}^m\frac{e^{i\vartheta_j/2}e^{x_j/2}}{e^{i\vartheta_j/2}e^{x_j/2}-e^{-i\vartheta_j/2}e^{-x_j/2}}=\prod_{j=1}^m\frac{1}{1-e^{-i\vartheta_j}e^{-x_j}}\]
Since the Chern roots of $N^*$ are $-x_j$, with $\gamma$-eigenvalues $e^{-i\vartheta_j}$,
\[ \ch_\gamma(\lambda_{-1}N^*)=\sum_k(-1)^k\ch_\gamma(\Lambda^kN^*)=\prod_{j=1}^m\bigl(1-e^{-i\vartheta_j}e^{-x_j}\bigr)\]
Together with Theorem~\ref{thm:Toeplitz character} this proves the proposition.

\end{proof}

Proposition~\ref{prop:Toeplitz integrand} shows that for Toeplitz operators Theorem~\ref{thm:main-intro} reads
\[ \Tr\bigl(\gamma\,|\,\Ind_GT_f\bigr)=\int_{M^\gamma}\ch_\gamma(f)\wedge\Td(M^\gamma)\wedge\frac{1}{\ch_\gamma(\lambda_{-1}N^*)}\]
For Toeplitz operators, Theorem~\ref{thm:main-intro} is therefore equivalent to the equivariant index formula which we derived in \cite[Theorem~2]{GvE26}.
This formula was proved by identifying 
the Szeg\"o projection with the Spin$^c$ Dirac operator of $M$ in equivariant $K$-homology, 
\[[D]=[S]\qquad \text{ in } K_1^G(M)\]
This  reduces the index of $T_f$ to the equivariant index of a twisted Spin$^c$ Dirac operator, and the Atiyah--Singer fixed-point formula gives the integrand.

\section{Proof of the main theorem}\label{sec:index-proof}

As in \cite{GvE20}, the proof is a reduction to Toeplitz operators.
If we fix $\gamma\in G$, then both sides of the index formula  depend only on the class $[\sigma_H(P)]\in K_1^G(\Symb)$.
It suffices to verify the formula for generators of $K_1^G(\Symb)$.
We show that $K_1^G(\Symb)$ is generated by two kinds of classes: symbols of multiplication operators and symbols of Toeplitz operators.
For multiplication operators both sides vanish.
For Toeplitz operators the formula is the equivariant index theorem of \cite{GvE26}, by Section~\ref{sec:Toeplitz operators}.

\subsection{Generators of $K_1^G(\Symb)$}\label{sec:K generation}

Let $W$ be a finite-dimensional unitary representation of $G$.
Multiplication by $g\in C^\infty(M,GL(W))^G$ is a $G$-equivariant order zero operator $\mathcal M_g$, and $\sigma_H(\mathcal M_g)=g$ via the inclusion $C^\infty(M)\otimes\End W\subset\Symb(M)\otimes\End W$ of Section~\ref{sec:Toeplitz set-up}.
In each fiber, $g$ is the distribution $g(x)\,\delta_e$.
We call $[g]\in K_1^G(\Symb)$ a \emph{multiplication class}.
For $f\in C^\infty(M,GL(W))^G$, the class
\[ [\sigma_H(\tilde T_f)]=[S_+\otimes f+(1-S_+)\otimes1_W]\in K_1^G(\Symb)\]
of Section~\ref{sec:Toeplitz set-up} is a \emph{Toeplitz class}.

\begin{proposition}\label{prop:K generation}
Every element of $K_1^G(\Symb)$ is of the form
\[ [g]+[\sigma_H(\tilde T_f)]\]
for a finite-dimensional unitary representation $W$ of $G$ and $f,g\in C^\infty(M,GL(W))^G$.
\end{proposition}

For trivial $G$ this is \cite[Proposition~8.1]{GvE20}.
We recall the argument and explain why it goes through equivariantly.

\begin{proof}[Sketch of proof]
By Section~\ref{sec:chi_gamma well defined} we may replace $\Symb(M)$ by its $C^*$-algebra completion $\overline{\Symb}(M)$.
Let $S^*M=S(H^*\oplus\RR)$ be the cosphere bundle, where $T^*M=H^*\oplus\RR\contact$, and let $\pi:S^*M\to M$ be the projection.
The sphere bundle $S^*H$ of $H^*$ is embedded in $S^*M$ as the equator $\{\xi(\Reeb)=0\}$.
For the metric of Section~\ref{sec:contact} the contact form is a unit covector, so $\pm\contact$ are sections of $S^*M$, through the two poles.
The open upper and lower hemispheres $\{\xi(\Reeb)>0\}$ and $\{\xi(\Reeb)<0\}$ are each identified with $H^*$.

The Fourier transform of a symbol, restricted to the hyperplane $\ft=0$, gives a homomorphism $\Symb(M)\to C(S^*H)$, which extends by continuity to $\overline{\Symb}(M)$.
Its kernel $I_H$ is the section algebra of a continuous field with fibers $C^*(H_x,\omega)\oplus C^*(H_x,-\omega)\cong\mathcal K\oplus\mathcal K$, the $\sigma_+$ and $\sigma_-$ components \cite[Section~8]{GvE20}.
We write $I_H=I_+\oplus I_-$, where $I_+$ is the summand that contains $S_+$, i.e.\ the $\sigma_+$ component.
In \cite[Section~8]{GvE20} (see \cite[Section~4.6, Proposition~4.6.1]{BvE14}), the extension
\[ 0\to I_H\to\overline{\Symb}(M)\to C(S^*H)\to0\]
is compared with
\[ 0\to C_0(H^*\sqcup H^*)\to C(S^*M)\to C(S^*H)\to0\]
where the two copies of $H^*$ are the open hemispheres.
The comparison is $C(M)$-linear and unital, it is the identity on the quotient $C(S^*H)$, and it maps each summand of $I_H$ to the functions on one hemisphere, by an isomorphism in $K$-theory (a Morita equivalence followed by a Thom isomorphism).
By the five lemma it gives an isomorphism $K_1(\overline{\Symb})\cong K^1(S^*M)$.

All of this is $G$-equivariant.
The comparison is constructed $U(n)$-equivariantly on the model $\Heis$ and transferred to $M$ through the unitary frame bundle $P_H$ of $\Hzo$ (Section~\ref{sec:curvature}).
The group $G$ acts on $P_H$ by $U(n)$-bundle maps.
Hence the comparison is $G$-equivariant, and the same argument gives isomorphisms in equivariant $K$-theory that make the diagrams
\[
\begin{CD}
K_1^G(C(M)) @>{j_*}>> K_1^G(\overline{\Symb})\\
@V{\cong}VV @VV{\cong}V\\
K^1_G(M) @>{\pi^*}>> K^1_G(S^*M)
\end{CD}
\qquad\qquad
\begin{CD}
K_1^G(I_+)\oplus K_1^G(I_-) @>>> K_1^G(\overline{\Symb})\\
@V{\cong}VV @VV{\cong}V\\
K^1_G(H^*)\oplus K^1_G(H^*) @>>> K^1_G(S^*M)
\end{CD}
\]
commute, where $j:C(M)\subset\overline{\Symb}(M)$ is the inclusion.

The proof is then a hemisphere argument.
Let $\iota\in\{\contact,-\contact\}$ be the pole section of the hemisphere that corresponds to $I_-$.
Since $G$ preserves $\contact$, $\iota$ is a $G$-equivariant section, and the closed hemisphere bundle that contains $\iota(M)$ retracts $G$-equivariantly onto $\iota(M)$.
For $x\in K^1_G(S^*M)$ let $g:=\iota^*x\in K^1_G(M)$.
Then $x-\pi^*g$ restricts to zero on $\iota(M)$, hence on this closed hemisphere bundle, and so it is in the image of $K^1_G(H^*)\to K^1_G(S^*M)$ from the other, open hemisphere.
By the first diagram $\pi^*g$ corresponds to the multiplication class $j_*g$.
By the second diagram $x-\pi^*g$ corresponds to an element of the image of $K_1^G(I_+)\to K_1^G(\overline{\Symb})$.

This image consists of Toeplitz classes.
The Szeg\"o symbol $S_+$ is a $G$-invariant projection in $I_+$, which in each fiber is the rank one vacuum projection in $C^*(H_x,\omega)\cong\mathcal K$ \cite[Section~3.6]{GvE26}.
Hence $S_+I_+S_+=C(M)S_+\cong C(M)$ is a full $G$-invariant corner of $I_+$, and the inclusion of the corner induces an isomorphism $K^1_G(M)\cong K_1^G(I_+)$, given by $[f]\mapsto[1_W+S_+\otimes(f-1_W)]$.
The image of this class in $K_1^G(\overline{\Symb})$ is $[S_+\otimes f+(1-S_+)\otimes1_W]=[\sigma_H(\tilde T_f)]$.

Finally, in $K^1_G(M)=K_1^G(C(M))$  every element is the class of a single invertible, since $[\sigma]-[\rho]=[\sigma\oplus\rho^{-1}]$.
A continuous $G$-equivariant $f:M\to GL(W)$ can be replaced by a smooth one.
Replacing $W$ by a direct sum $W_1\oplus W_2$ and $f$, $g$ by $f\oplus1$, $g\oplus1$, we may take the same $W$ for $f$ and $g$.

\end{proof}

\subsection{Conclusion of the proof}\label{sec:conclusion}

\begin{proof}[Proof of Theorem~\ref{thm:main-intro}]
It suffices to prove the formula for the classes $[g]$ and $[\sigma_H(\tilde T_f)]$ for $P=\mathcal M_g$ and for $P=\tilde T_f$.

For $P=\mathcal M_g$ both sides vanish, by Lemma~\ref{lem:multiplication classes} below.

For $P=\tilde T_f$  the right-hand side of the formula in Theorem~\ref{thm:main-intro}
is equivalent to the equivariant index formula for Toeplitz operators.

\end{proof}

\begin{lemma}\label{lem:multiplication classes}
For any $g\in C^\infty(M,GL(W))^G$ 
we have  $\chi_\gamma(g)=0$ as a differential form on $M^\gamma$.
\end{lemma}
\begin{proof}
Consider the factors in the terms of Definition~\ref{def:chi_gamma}, with $\sigma=g$.
In each fiber, $\sigma=g(x)\delta_e$ and $\sigma^{-1}=g(x)^{-1}\delta_e$.
In a local frame, $\beta\cdot\delta_e=0$, because $\Sp(V)$ acts on $\Heis$ by group automorphisms that fix $e$.
Hence $\con\sigma=d\sigma+\beta\cdot\sigma=dg\,\delta_e$, and likewise $\con(\sigma^{-1})=d(g^{-1})\,\delta_e$.
The curvature $\curv=\nu(\theta)$ takes values in $\nu(\Lu(\Hzo_x))\subset\Tsf^{-1}\mcU^2(\heis)$ (Section~\ref{sec:enlarged-symbol-algebra}).
So at each point of $M^\gamma$ every term
\[ \sigma^{-1}\curv^{i_0}\con(\sigma)\curv^{i_1}\con(\sigma^{-1})\cdots\con(\sigma)\curv^{i_{2l+1}}\]
is a form whose coefficients are elements of $\Alg_H(\Heis_x)\otimes\End W$ supported at the identity.

On such elements $\orbtr{c}$ vanishes for $c\ne0$,
since  the orbit $\Gamma_c=\{(z,c+q(z))\mid z\in N\}$ does not contain $e=(0,0)$.
Hence $\orbtr{\pm1}$, and with them $\Tr_\gamma$ vanish on every term.

\end{proof}

\appendix
\section{Conventions}\label{sec:conventions}

\subsection*{Signs and constants}

Table~\ref{tab:conventions} collects the conventions on which the signs and constants of the index formula depend.

\begin{table}[ht]
\centering
\small
\renewcommand{\arraystretch}{1.7}
\begin{tabular}{@{}>{\raggedright}p{0.30\textwidth}>{\raggedright\arraybackslash}p{0.62\textwidth}@{}}
\toprule
Symplectic form & $\omega=d\contact|_H$\\
Metric &  $g(v,v)=\omega(v,Jv)>0$, $\ang{v,w}=g(v,w)-i\omega(v,w)$ \\
Group law & $(z,t)(w,s)=(z+w,t+s-\frac12\omega(z,w))$\\
Measure & $\dd h=\dd z\,\dd t$, $\dd z=|\omega^n|/n!$ \\
Orientations & $M$: $\contact\wedge(d\contact)^n>0$; components of $M^\gamma$: $\contact\wedge(d\contact)^r>0$ \\
Homogeneity & $\delta_a(z,t)=(az,a^2t)$, $Q=2n+2$; $\cR^j$: degree $-Q-j$ \\
Central element & $\Tsf=-\partial_t\delta_e$, $\Tsf*f=-\partial_tf$, $\Tsf\in\cR^2$ \\
Fourier transform & $\mcF f(\xi,\ft)=\int f(z,t)\,e^{-i(\xi(z)+\ft t)}\,\dd z\,\dd t$; $\mcF\Tsf=-i\ft$ \\
Weyl pairs & $\sigma_+=\mcF f|_{\ft=-1}$, $\sigma_-=\mcF f|_{\ft=+1}$; $\sigma_H(\Reeb)=(i,-i)$ \\
Szeg\"{o} symbol & Weyl pair  $S_+=(s_0,0)$, $s_0(\xi)=2^ne^{-|\xi|^2}$\\
Szeg\"{o} kernel & $\frac{n!}{(2\pi)^{n+1}}\bigl(\frac{|z|^2}{4}+it\bigr)^{-n-1}$ of the identity\\
Eigenvalues of $\gamma_N$ & $e^{i\vartheta_j}$, $0<\vartheta_j<2\pi$; $\sqrt{\det\gamma_N}=\prod_je^{i\vartheta_j/2}$ \\
Twisted orbits & $\Gamma_c=\{(z,c+q(z))\}$, $q(z)=\ang{Cz,z}$, $C=\frac i4(\gamma_N+1)(\gamma_N-1)^{-1}$ \\
Orbital traces & $\orbtr{c}(\Tsf^{-j}f)=\frac{r!\,c^j}{(r+j)!}\int_Nf(z,c+q(z))\,\dd z$ \\
Equivariant trace & $\Tr_\gamma=\kappa_\gamma\frac{(2\pi i)^{r+1}}{2(r!)}\bigl(\orbtr{+1}-(-1)^r\orbtr{-1}\bigr)$, $\kappa_\gamma=\det\frac{\gamma_N^{1/2}}{1-\gamma_N}$ \\
Values & $\Tr_\gamma(S_+)=\sqrt{\det\gamma_N}$\\
Vacuum curvature & $\nu(A)S_+=\frac12\tr(A)S_+$ \\
Character & coefficients $\bigl(\frac{-1}{2\pi i}\bigr)^{I+l+1}\frac{l!}{(I+2l+1)!}$, i.e.\ $u=\frac{1}{2\pi i}$ \\
Chern--Weil & $\ch_\gamma(E)=\tr\bigl(\gamma\exp(\frac{i}{2\pi}\Omega^E)\bigr)$, $c_1(E)=[\frac{i}{2\pi}\tr\Omega^E]$ \\
Fixed-point denominator & $\Delta_\gamma(N)=\prod_j\bigl(e^{i\vartheta_j/2}e^{x_j/2}-e^{-i\vartheta_j/2}e^{-x_j/2}\bigr)$ \\
\bottomrule
\vskip 6pt
\end{tabular}
\caption{Signs and constants. Here $r=\dim_\CC\ker(1-\gamma)$, $\gamma_N$ is the action of $\gamma$ on the normal part $N$, and the $x_j$ are the Chern roots of $N$.}
\label{tab:conventions}
\end{table}

\subsection*{Comparison with \cite{GvE20}, \cite{GvE25} and \cite{GvE26}}

The symplectic form $\omega=d\contact|_H$ agrees with \cite{GvE26}.
In \cite{GvE20} the symplectic form is $-d\contact$; the complex structure $J$, the metric $g$, the bundle $\Hzo$ and the bracket of $\heis$ are the same as ours, and only the name $\omega$ has changed sign.
In \cite{GvE25} the group law is $(x,t)(y,s)=(x+y,t+s+\frac12\omega(x,y))$, so the symplectic form there is also $-\omega$.

Our Weyl pair labels are those of \cite{GvE20} and \cite[Section~4.2]{GvE26}, with $\sigma_H(\Reeb)=(i,-i)$.
In \cite{GvE25}, $\sigma_\pm$ denotes the restriction to the slice $\ft=\pm1$, i.e.\ our $\sigma_\mp$.
Together with the opposite sign of $\omega$, this accounts for the sign by which the formula for $\tau$ in \cite{GvE25} differs from the one in Section~\ref{sec:Tr_gamma}.

\subsection*{Equivariant Chern characters}

We use the conventions of \cite[Section~5.7]{GvE26}.
Let $\gamma\in G$.
For a $G$-equivariant Hermitian vector bundle $E$ on $M$ with a $G$-invariant unitary connection with curvature $\Omega^E$,
\[ \ch_\gamma(E)=\tr\Bigl(\gamma\exp\bigl(\tfrac{i}{2\pi}\Omega^E\bigr)\Bigr)\qquad c_1(E)=\Bigl[\tfrac{i}{2\pi}\tr(\Omega^E)\Bigr]\]
on $M^\gamma$.
For a finite-dimensional unitary representation $W$ of $G$ and a smooth $G$-equivariant map $f:M\to GL(W)$, the odd equivariant Chern character is represented by
\[ \ch_\gamma(f)=-\sum_{l\ge0}\frac{l!}{(2l+1)!\,(2\pi i)^{l+1}}\,\tr_W\bigl((f^{-1}df)^{2l+1}\gamma\bigr)\in\Omega^{\mathrm{odd}}(M^\gamma)\]
For $\gamma=e$ this is the odd Chern character of \cite{GvE20}, denoted $\Ch(f)$ there.

\bibliographystyle{amsplain}
\bibliography{MyBibfile2}
\end{document}

%% file: macros.tex
\usepackage{mathrsfs}   

\providecommand{\CC}{{\mathbb{C}}}
\providecommand{\RR}{{\mathbb{R}}}

\providecommand{\ZZ}{{\mathbb{Z}}}

\providecommand{\mcE}{{\mathcal{E}}}
\providecommand{\mcF}{{\mathcal{F}}}   

\providecommand{\mcU}{{\mathcal{U}}}   

\newcommand{\contact}{\alpha}                    
\newcommand{\Reeb}{T}                            
\newcommand{\Hzo}{H^{1,0}}                       

\providecommand{\Heis}{\mathbb{H}}               
\providecommand{\heis}{\mathfrak{h}}             
\providecommand{\Lg}{\mathfrak{g}}               
\providecommand{\Lu}{\mathfrak{u}}
\providecommand{\Sp}{\mathrm{Sp}}
\providecommand{\spl}{\mathfrak{sp}}
\newcommand{\Tsf}{\mathsf{T}}                    
\newcommand{\ft}{s}                              

\providecommand{\Symb}{{\mathscr{S}_H}}          

\providecommand{\Alg}{{\mathscr{A}}}             
\providecommand{\cR}{\mathcal{R}}                
\providecommand{\Trh}{{\overline{\mathrm{Tr}}}}  

\newcommand{\con}{\boldsymbol{\nabla}}           
\newcommand{\curv}{\boldsymbol{\Theta}}          

\newcommand{\orbtr}[1]{\tau_{#1}}                
\newcommand{\orbint}[1]{I_{#1}}                  

\newcommand{\Ahat}{\widehat{A}}                 

\providecommand{\ang}[1]{\langle #1 \rangle}

\providecommand{\dd}{\,d}                        

\providecommand{\End}{\operatorname{End}}
\providecommand{\Tr}{\operatorname{Tr}}
\providecommand{\tr}{\operatorname{tr}}
\providecommand{\Ch}{\operatorname{Ch}}          
\providecommand{\ch}{\operatorname{ch}}          
\providecommand{\Td}{\operatorname{Td}}          
\providecommand{\ind}{\operatorname{Index}}      
\providecommand{\Ind}{\operatorname{Index}}      

\providecommand{\supp}{\operatorname{supp}}
